\RequirePackage{fix-cm}
\documentclass{article}

\usepackage{arxiv}

\usepackage[utf8]{inputenc} 
\usepackage[T1]{fontenc}    
\usepackage[colorlinks=true,urlcolor=black,citecolor=black,linkcolor=black]{hyperref}
\usepackage{url}            
\usepackage{booktabs}       
\usepackage{nicefrac}       

\usepackage{graphicx}
\usepackage{doi}

\usepackage{color}
\usepackage{amssymb,amsfonts,amsmath,amsthm,amsbsy} 
\usepackage{tabularx} 
\usepackage{mathtools,mathrsfs}
\usepackage{mhchem}
\usepackage{enumerate}
\usepackage{bm}
\usepackage{bbm}
\usepackage{multirow}
\usepackage{threeparttable}
\usepackage{float}

\usepackage[table,dvipsnames,svgnames]{xcolor}
\definecolor{darkgreen}{rgb}{0.0, 0.5, 0.0}

\usepackage[bottom]{footmisc}
\usepackage[numbers,sort&compress]{natbib}

\usepackage{indentfirst}
\usepackage{etoolbox,siunitx}
\renewcommand{\thefootnote}{\fnsymbol{footnote}}
\usepackage{afterpage}
\robustify\bfseries
\usepackage{adjustbox}
\usepackage{changepage}
\usepackage{rotating}
\usepackage{hvfloat}
\usepackage{soul}
\usepackage{pifont}

\usepackage{placeins}
\usepackage{todonotes}

\usepackage{subcaption}
\usepackage{caption}

\usepackage{algorithm}
\usepackage{algpseudocode}
\usepackage[italicComments=true,indLines=true,noEnd=false]{algpseudocodex}
\usepackage[english]{babel}
\hypersetup{colorlinks=true, citecolor=blue}
\usepackage{stackrel}
\usepackage[nopatch=eqnum]{microtype}      

\usepackage[framemethod=TikZ]{mdframed}

\mdfdefinestyle{stepbox}{
    roundcorner=7pt,
    linewidth=1pt,
    innertopmargin=1pt,
    innerbottommargin=2pt,
    innerleftmargin=5pt,
    innerrightmargin=5pt,
    skipabove=1pt,
    skipbelow=1pt
}

\mdfdefinestyle{gray}{style=stepbox, linecolor=black!20, backgroundcolor=gray!5}
\mdfdefinestyle{dmd_blue}{style=stepbox, linecolor=cyan!60!blue!40, backgroundcolor=cyan!7}
\mdfdefinestyle{sindy_purple}{style=stepbox, linecolor=violet!40, backgroundcolor=violet!7}
\mdfdefinestyle{jacobian_orange}{style=stepbox, linecolor=orange!50, backgroundcolor=orange!10}
\mdfdefinestyle{csp_teal}{style=stepbox, linecolor=teal!40, backgroundcolor=teal!10}
\mdfdefinestyle{sim_green}{style=stepbox, linecolor=green!50!teal!40, backgroundcolor=green!7}
\mdfdefinestyle{rom_magenta}{style=stepbox, linecolor=magenta!40, backgroundcolor=magenta!5}

\definecolor{abstract_background}{RGB}{235,235,235}%
\usepackage[title]{appendix}
\usepackage[capitalise,nameinlink]{cleveref}

\title{Data-driven framework for invariant manifold identification in multiscale systems}

\author{
\textbf{\href{https://orcid.org/0000-0002-4788-9877}{\includegraphics[scale=0.06]{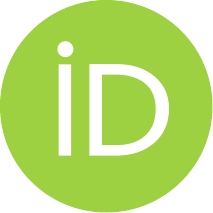}\hspace{1mm}Dimitris M. Manias}}\textcolor{teal}{$^{1,\ddagger}$}\thanks{Corresponding Author: dimitris.manias@ku.ac.ae}~,
\textbf{\href{https://orcid.org/0009-0005-6443-2894}{\includegraphics[scale=0.06]{orcid.pdf}\hspace{1mm}Ismaila Muhammed}}\textcolor{teal}{$^{1,\ddagger}$}~,
\textbf{\href{https://orcid.org/0000-0002-9674-9491}{\includegraphics[scale=0.06]{orcid.pdf}\hspace{1mm}Dimitris A. Goussis}}\textcolor{teal}{$^{2}$}~,
\textbf{\href{https://orcid.org/0000-0002-1270-7885}{\includegraphics[scale=0.06]{orcid.pdf}\hspace{1mm}Haralampos Hatzikirou}}\textcolor{teal}{$^{3,*}$}\thanks{Corresponding Author: h.hatzikirou@euc.ac.cy}
{}\\[5pt]
\textcolor{teal}{$^{(1)}$} Department of Mathematics, \emph{Khalifa University of Science and Technology}, \\ Abu Dhabi 127788, UAE \\
\textcolor{teal}{$^{(2)}$} Department of Mechanical Engineering, \emph{Khalifa University of Science and Technology}, \\ Abu Dhabi 127788, UAE \\
\textcolor{teal}{$^{(3)}$} Life and Health Sciences, \emph{European University of Cyprus}, Frankfurt, Germany
}

\renewcommand{\shorttitle}{Data-driven framework for invariant manifold identification in multiscale systems}

\hypersetup{
pdftitle={Data-driven framework for invariant manifold identification in multiscale systems},
pdfsubject={math.DS, q-bio.QM},
pdfauthor={Dimitris M. Manias, Ismaila Muhammed, Dimitris A. Goussis, Haralampos Hatzikirou},
pdfkeywords={Slow invariant manifold, Computational Singular Perturbation, SINDy, Krylov-Arnoldi, Mori-Zwanzig, multiscale dynamical systems},
}

\begin{document}
\maketitle

\renewcommand{\thefootnote}{\fnsymbol{footnote}}
\footnotetext[3]{These authors contributed equally as co-first authors.}
\renewcommand{\thefootnote}{\arabic{footnote}}

\begin{abstract}
\colorbox{abstract_background}{\begin{minipage}{1\linewidth}
Multiscale dynamical systems often exhibit fast transient behavior followed by slow evolution on a low-dimensional attracting manifold, on which their long-term dynamics unfold. Identifying this structure directly from data remains challenging when the governing equations, Jacobians, and fast-slow variable partitions are unknown. In this work, we propose a fully data-driven framework for timescale decomposition, slow invariant manifold (SIM) approximation, and reduced-order model identification. A sparse-regression surrogate of the vector field is reconstructed from the data, and its Jacobian, estimated with a matrix-free Krylov method, is analysed with Computational Singular Perturbation diagnostics. The resulting fast-slow partition guides the learning of explicit SIM algebraic representations through regression, and the reduced model follows by restricting the dynamics to it. For a linear stochastic system, this regression agrees with the Mori-Zwanzig projection at leading order in the timescale ratio. We demonstrate the validity of the framework on an illustrative linear system, the Michaelis-Menten model, a biofilm growth model, and a five-dimensional target-mediated drug disposition (TMDD) model. The method accurately identifies fast variables, reconstructs slow invariant manifolds, and generates accurate reduced models. Furthermore, we map the performance boundaries of the framework under sparse sampling and measurement noise, showing robust recovery for low-dimensional systems and highlighting data-requirement limits in higher dimensions. These results demonstrate that interpretable multiscale reduction can be achieved directly from time-series data, without requiring explicit prior knowledge of the underlying fast/slow timescale separation.
\end{minipage}
}
\end{abstract}
\par\noindent\rule{\textwidth}{0.5pt}

\keywords{Slow invariant manifold \and Computational Singular Perturbation \and SINDy \and Krylov-Arnoldi \and Mori-Zwanzig \and multiscale dynamical systems}


\section{Introduction}

Biological systems exhibit complex behavior governed by interacting components and biochemical processes \cite{akbari2024data}.~Detailed kinetic models formulated from enzymatic rate laws can, in principle, describe these interactions with high fidelity; however, their application to large-scale biochemical systems is often limited by numerical stiffness, high computational cost, and reduced interpretability \cite{karr2012whole,macklin2020simultaneous}.~Despite their complexity, such systems frequently exhibit relatively simple multiscale behavior arising from underlying timescale hierarchies, ranging from nanoseconds (\(10^{-9}\) s) at the molecular scale to years (\(10^{8}\) s) at the organism level \cite{heinrich1982dynamics,qu2011multi}.~Consequently, timescale decomposition has become an important framework for the analysis and development of reduced-order models that retain the essential dynamics governing the long-time behavior of these systems.

Existing approaches to timescale decomposition can generally be classified as either top-down or bottom-up \cite{akbari2024data, qu2011multi}.~Top-down methods infer correlated groups of species or states, often described as concentration pools, directly from time-series observations; however, they typically lack a direct mechanistic connection to the underlying stoichiometric or dissipative structure of the network \cite{arkin1995statistical,arkin1997test}.~Bottom-up methods, in contrast, extract timescale information from eigenvalue analysis of the system Jacobian and therefore retain a closer connection to the governing equations; however, they require explicit knowledge of reaction-rate laws and may become inaccurate in highly nonlinear systems whose Jacobian spectra vary significantly over time \cite{kauffman2002description,jamshidi2008top}.

More broadly, dimensionality reduction plays a central role in timescale decomposition \cite{radulescu2008robust}. Trajectory-based reduction techniques can capture rich time-dependent behavior without imposing an explicit fast–slow decomposition, but they are often computationally demanding, sensitive to measurement noise, and reliant on large volumes of data \cite{reinhardt2008approximation,schmid2022dynamic}. Singular perturbation methods provide a rigorous mathematical framework for separating fast and slow dynamics; however, their practical application remains challenging due to their reliance on explicit reaction-rate expressions and kinetic parameters \cite{lam1994csp,chiavazzo2007comparison}. 

In many multi-scale dynamical systems, trajectories exhibit a characteristic two-stage evolution: an initial rapid relaxation toward an attracting low-dimensional manifold, followed by slower evolution along this manifold \cite{fenichel1979geometric, jones2006geometric, gear2005projecting}.~This emergent structure is commonly referred to as the \emph{slow invariant manifold} (SIM), and its accurate approximation is central to both mechanistic understanding and the construction of parsimonious reduced-order models \cite{kuehn2015multiple,zagaris2004analysis}.

Classical approaches to SIM approximation are rooted in singular perturbation theory, geometric singular perturbation theory, and numerical solutions of the invariance equation \cite{fenichel1979geometric,jones2006geometric,gorban2004constructive,verhulst2005methods}.~These methods can only be applied for the analysis of relatively simple models, which can be cast in the appropriate non-dimensional form, so that the fast–slow decomposition and variable partition are known \textit{a priori}.~For complex systems, in which such information is not easily available, the Computational Singular Perturbation (CSP) method and its variances, Intrinsic Low-Dimensional Manifolds (ILDM) and Tangential Stretch Rate (TSR), provide algorithmically approximations of the SIM and of the reduced system that governs the flow along the SIM,  by exploiting the local tangent-space structure \cite{goussis2012quasi,maas1992simplifying,valorani2017tangential}.

The use of model reduction methodologies is most promising in handling the increasingly available time-series data.~Classical reduced-order modeling approaches, including projection-based and subspace identification methods, represent the reduced dynamics in transformed coordinates, thereby sacrificing direct physical interpretability \cite{benner2015survey,overschee1996subspace, holmes2012turbulence,sirovich1987turbulence, moore2003principal, mullis1976synthesis,rowley2009spectral, schmid2010dynamic, juang1985eigensystem}.~The recent machine learning–based approaches have demonstrated strong predictive capabilities, and often operate as black-box models without consider the time-scale separation and the SIM  \cite{gallos2024data,chen2021physics, williams2015data, wan2018data, demo2023deeponet, vlachas2022multiscale, vijayarangan2024data, gruber2022comparison}. Recent physics-informed approaches either assume an explicit singularly perturbed form with known fast and slow variables or, when they infer the decomposition, require the governing equations \cite{patsatzis2024slow, serino2025fast, patsatzis2024physics}.

At this stage, the existing approaches for model identification require the availability of analytical models, assume known timescale separation, or fail to preserve interpretability in high-dimensional settings.~Therefore, the identification of multi-scale structure and reduced model directly from existing data sets, without prior knowledge of governing equations or fast/slow dynamics, remains a major challenge.

In this work, we propose a unified data-driven framework for system identification in a multi-scale context, that overcomes the existing limitations by integrating sparse system identification, matrix-free Jacobian estimation, and timescale decomposition.~Specifically, we combine Sparse Identification of Nonlinear Dynamics (SINDy) with a Krylov subspace–based Arnoldi method to estimate the local Jacobian directly from data, thus enabling the extraction of timescale information without explicit knowledge of the governing equations.~We then employ CSP diagnostics to identify fast and slow dynamics, allowing for the detection of fast variables and the approximation of the SIM.~The SIM is subsequently learned using regression techniques, and reduced-order models are identified using Weak SINDy restricted to the slow manifold.~This work is an extension of our previous work for the algorithmic identification of valid models locally from biological data \cite{muhammed2025data}.

The proposed framework provides an interpretable, fully data-driven pipeline for multi-scale model reduction that preserves physical meaning, while capturing the essential system dynamics.~Its effectiveness is demonstrated on several benchmark systems, including a linear 2-dimensional model, the 2-dimensional Michaelis--Menten (MM) model, a 3-dimensional biofilm growth dynamics model, and a 5-dimensional target-mediated drug disposition (TMDD) model.~It is shown that the framework successfully identifies fast variables, reconstructs the slow manifold, and yields accurate reduced models.

The paper is organized as follows.~Section 2 introduces the proposed framework for timescale decomposition and reduced modeling.~Section 3 illustrates its application to the Michaelis--Menten, biofilm, and TMDD models.~Section 4 concludes the paper and outlines directions for future research.

\section{Proposed framework}

We propose a fully data-driven framework for the identification of SIMs and reduced models in nonlinear multi-scale dynamical systems directly from time-series data.~The framework integrates (i) temporal data augmentation, (ii) sparse vector-field reconstruction, (iii) matrix-free Jacobian estimation, and (iv) CSP diagnostics.

The framework is organized into two main stages.~First, sparse or heterogeneous time-series data are transformed into a smooth surrogate dynamical system through DMD and sparse regression.~This system is then employed for matrix-free Jacobian estimation.~Second, the estimated Jacobian is analyzed using CSP diagnostics to identify time-scale separation, fast variables, and asymptotic constraints.~These identifications enable regression-based approximation of the slow invariant manifold and identification of reduced model.

\begin{figure}[h]
\centering
\includegraphics[width=\linewidth]{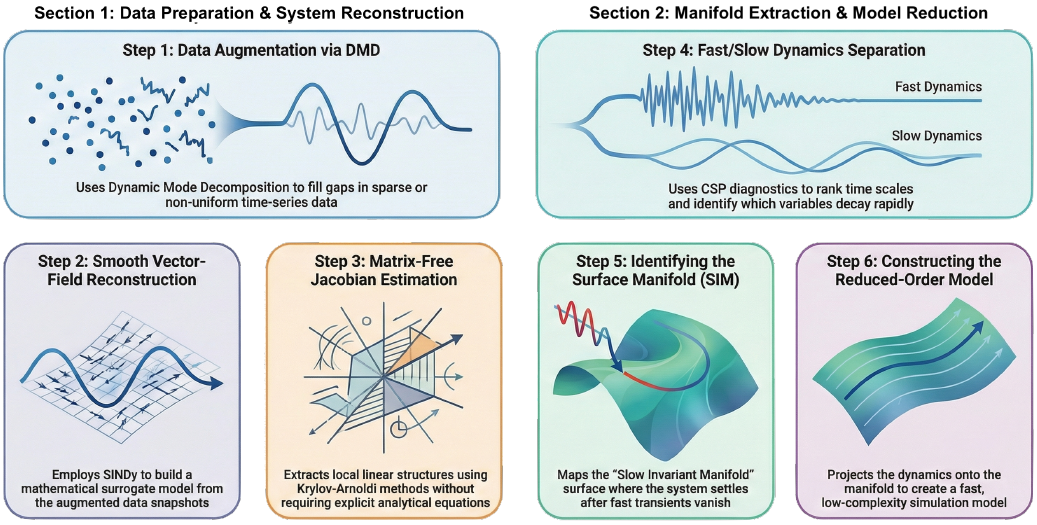}
\caption{Proposed data-driven framework for slow invariant manifold identification and reduced-order modeling.~Sparse or heterogeneous time-series data are first augmented using DMD (Step 1), followed by sparse vector-field reconstruction (Step 2) and matrix-free Jacobian estimation through Krylov--Arnoldi methods  (Step 3).~The estimated Jacobian  is subsequently analyzed using CSP diagnostics to identify fast and slow dynamics (Step 4), approximate SIMs  (Step 5), and construct reduced models approximating the flow on the learned manifold  (Step 6).}
\label{fig:framework}
\end{figure}
\FloatBarrier

\subsection{Data-driven system reconstruction}

We consider physical processes that can be approximated by  nonlinear multi-scale dynamical systems of the form:
\[
\dot{\mathbf{y}}(t)=\mathbf{g}(\mathbf{y}(t)),
\qquad
\mathbf{y}(t)\in\mathbb{R}^{N},
\]
where the vector field $\mathbf{g}$ is unknown or only partially known, and only time-resolved observations of the state $\mathbf{y}(t)$ are available.~Let $\{\mathbf{y}(t_k)\}_{k=1}^{T_0}$ denote the original set of observed snapshots, where the sampling times $\{t_k\}$ may be sparse or nonuniform.~To improve temporal resolution and consistency, DMD is employed to approximate the evolution of the system through a linear operator:
\[
\mathbf{y}(t_{k+1})
\approx
\mathbf{A}\,\mathbf{y}(t_k),
\]
where $\mathbf{A}\in\mathbb{R}^{N\times N}$ is inferred directly from the data. The spectral decomposition of $\mathbf{A}$ provides dynamically consistent reconstructions at arbitrary intermediate times, yielding an augmented trajectory
\[
\{\widetilde{\mathbf{y}}(t_k)\}_{k=1}^{T},
\qquad
T\gg T_0,
\]
defined on a prescribed reference grid. In the presence of heterogeneous sampling across trajectories, DMD is applied independently to each realization prior to aggregation on a common time grid. Although classical DMD provides an efficient linear approximation of the underlying dynamics, its accuracy may deteriorate for strongly nonlinear systems because it assumes that the state evolves under a linear operator. To improve the representation of nonlinear dynamics, Extended Dynamic Mode Decomposition (EDMD) augments the original state variables with a dictionary of nonlinear observables, thereby approximating a linear Koopman operator in the lifted feature space \cite{williams2015data}. More recently, Li et al. \cite{li2017extended} proposed Extended Dynamic Mode Decomposition with Dictionary Learning (EDMD-DL), in which the dictionary is learned directly from data using neural networks rather than being specified a priori. By adaptively constructing nonlinear observables, EDMD-DL improves the approximation of nonlinear dynamics while requiring substantially smaller dictionaries than conventional EDMD. These developments illustrate the growing trend toward data-adaptive representations of nonlinear systems.

For notational simplicity, the augmented dataset is subsequently represented by the data matrix:
\[
\mathbf{Y}
=
\begin{bmatrix}
\mathbf{y}(t_1),
\mathbf{y}(t_2),
\cdots ,
\mathbf{y}(t_T)
\end{bmatrix}
\in
\mathbb{R}^{N\times T},
\]
where each column corresponds to the system state at a sampling instant and each row corresponds to a specific state variable.~A smooth surrogate vector field is then reconstructed from the augmented trajectories using Sparse Identification of Nonlinear Dynamics (SINDy). Specifically, the dynamics are approximated as
\begin{equation}
\label{construct_g}
\hat{\mathbf{g}}(\mathbf{y})
=
\boldsymbol{\Theta}(\mathbf{y})\,\boldsymbol{\Xi},
\end{equation}
where $\boldsymbol{\Theta}(\mathbf{y})$ denotes a library of candidate nonlinear functions and $\boldsymbol{\Xi}$ is a sparse coefficient matrix identified through regularized sparse regression methods such as SR3 or sequential thresholded least squares \cite{brunton_discovering_2016}.~The resulting surrogate system:
\begin{equation}
\label{eq:reconstructed_vector_field}
\dot{\mathbf{y}}
=
\hat{\mathbf{g}}(\mathbf{y})
\end{equation}
provides a smooth and differentiable approximation of the dynamics suitable for local linear analysis.

Since analytic expressions for the Jacobian are not assumed to be available, local Jacobian information is extracted through a matrix-free Krylov--Arnoldi procedure~\cite{knoll2004jacobian}.~For a given state $\mathbf{y}\in\mathbb{R}^{N}$ and direction $\mathbf{v}\in\mathbb{R}^{N}$, Jacobian-vector products are approximated through finite differences:
\[
\mathbf{J}\mathbf{v}
\approx
\frac{
\hat{\mathbf{g}}(\mathbf{y}+\varepsilon\mathbf{v})
-
\hat{\mathbf{g}}(\mathbf{y})
}{
\varepsilon
},
\]
where $\varepsilon>0$ is a small perturbation parameter.~Using this operator, the Arnoldi iteration constructs an $r$-dimensional Krylov subspace:
\[
\mathcal{K}_{r}(\mathbf{J},\mathbf{v}_0)
=
\mathrm{span}
\{
\mathbf{v}_0,
\mathbf{J}\mathbf{v}_0,
\dots,
\mathbf{J}^{r-1}\mathbf{v}_0
\},
\]
together with an orthonormal basis $\mathbf{V}_{r}$ and upper Hessenberg matrix $\mathbf{H}_{r}$.~This yields the low-rank Jacobian approximation:
\begin{equation}
\label{eq:low_rank_jac}
\mathbf{J}_{r}
\approx
\mathbf{V}_{r}
\mathbf{H}_{r}
\mathbf{V}_{r}^\top,
\end{equation}
whose dominant eigenvalues approximate the eigenvalues of the true Jacobian.~Applying this procedure along the trajectory produces a sequence of local Jacobian approximations
\[
\{
\mathbf{J}_{r}(\mathbf{y}(t_k))
\}_{k=1}^{T},
\]
which form the basis for the subsequent multi-scale analysis.

\subsection{CSP-inspired multiscale diagnostics}

In the CSP framework, a multiscale vector field admits the modal decomposition:
\begin{equation}
\label{eq:CSP2}
\mathbf{g}(\mathbf{y})
=
\sum_{n=1}^{N}
\boldsymbol{\alpha}_n
f^n(\mathbf{y}),
\qquad
f^n(\mathbf{y})
=
\boldsymbol{\beta}^n
\mathbf{g}(\mathbf{y}),
\end{equation}
where $\{\boldsymbol{\alpha}_n\}_{n=1}^{N}$ and $\{\boldsymbol{\beta}^n\}_{n=1}^{N}$ denote biorthogonal right and left CSP basis vectors satisfying:
\[
\boldsymbol{\beta}^i\cdot\boldsymbol{\alpha}_j
=
\delta_j^i.
\]
Each term $\boldsymbol{\alpha}_n f^n$ is referred to as a CSP mode and consists of the direction defined by the basis vector $\boldsymbol{\alpha}_n$ together with the modal amplitude $f^n$, which quantifies the projection of the vector field $\mathbf{g}(\mathbf{y})$ onto $\boldsymbol{\alpha}_n$ \cite{lam1989understanding,hadjinicolaou1998asymptotic}.~Each mode is associated with a time scale $\tau_n$, ordered such that:
\[
\tau_1 < \tau_2 < \dots < \tau_N.
\]
The modal amplitudes quantify the contribution of the corresponding modes to the system evolution, whereas the associated time scales characterize the time frame of their activity.~These time scales may be either explosive or dissipative, indicating whether the corresponding modes locally drive the system away from or toward equilibrium, respectively \cite{diamantis2015,khalil2021no}.

In CSP, the basis vectors $\boldsymbol{\alpha}_n$ and $\boldsymbol{\beta}^n$ are obtained through an iterative refinement procedure that yields higher-order asymptotic corrections \cite{zagaris2004fast,kaper2015geometry}.~In the present framework, the right and left eigenvectors of the estimated Jacobian,
\[
\{
\mathbf{a}_n(t_k),
\mathbf{b}^n(t_k)
\}_{n=1}^{N},
\]
are interpreted as leading-order CSP basis approximations \cite{lam1989understanding,kuehn2015multiple}.

Let $\{\lambda_n(t_k)\}_{n=1}^{N}$ denote the eigenvalues of the estimated Jacobian at time $t_k$. The associated developing time scales are defined as
\begin{equation}
\label{eq:timescales}
\tau_n(t_k)
=
|\lambda_n(t_k)|^{-1},
\end{equation}
providing a local fast-to-slow ordering of the dynamics:
\[
\tau_1(t_k)
<
\tau_2(t_k)
<
\dots
<
\tau_N(t_k).
\]
Modes with $\Re(\lambda_n(t_k))<0$ are classified as dissipative, whereas modes with $\Re(\lambda_n(t_k))>0$ are locally explosive.

To identify the variables driven by a specific mode, the CSP pointer is employed:
\begin{equation}
\label{eq:Po}
Po^n
=
\operatorname{diag}
\left[
\mathbf{a}_n(t_k),
\mathbf{b}^n(t_k)
\right]
=
\left[
b_1^n a_n^1,
\dots,
b_N^n a_n^N
\right],
\end{equation}
where $a_n^i$ and $b_i^n$ denote the $i$-th components of the right and left eigenvectors, respectively, and $\sum_{i=1}^N b^n_i a^i_n = 1$ \cite{goussis1992}.~The magnitude of the $i$-th component of $Po^n$
assesses the way the $i$-th state variable will adjust when the system is perturbed along the direction $\boldsymbol{a}_n$ \cite{manias2025physical}.~Large pointer values indicate strong influence of the mode to the variable $y_i$.~Variables that consistently are pointed by the fast modes are identified as fast variables \cite{manias2016mechanism,manias2022effect,rabbani2022dominant}

Combined analysis of modal amplitudes, time scales, and pointer values of the fast modes provides qualitative information regarding the structure of fast relaxation.~The acquired information can be used to assess the possible validity of  quasi-steady-state (QSS) or partial-equilibrium (PE) behavior \cite{goussis2012quasi}.

\subsection{Slow invariant manifold identification and reduced-order modeling}

Assuming that the system exhibits $M$ fast dissipative time scales satisfying:
\[
\tau_1(t_k)
<
\dots
<
\tau_M(t_k)
\ll
\tau_{M+1}(t_k)
<
\dots
<
\tau_N(t_k),
\]
the CSP decomposition of Eq.~\eqref{eq:CSP2} allows for the separation of the reconstructed dynamics into fast and slow components:
\begin{equation}
\label{eq:decomposition}
\hat{\mathbf{g}}(\mathbf{y})
=
\hat{\mathbf{g}}_f(\mathbf{y})
+
\hat{\mathbf{g}}_s(\mathbf{y})
=
\sum_{m=1}^{M}
\mathbf{a}_m
f^m(\mathbf{y})
+
\sum_{n=M+1}^{N}
\mathbf{a}_n
f^n(\mathbf{y}).
\end{equation}
After the fast dissipative time scales become exhausted, the corresponding fast modal amplitudes satisfy:
\begin{equation}
\label{eq:CSP_sim}
f^m(\mathbf{y})
\approx
0,
\qquad
m=1,\dots,M,
\end{equation}
yielding the slow approximation:
\begin{equation}
\label{eq:CSP_slow}
\hat{\mathbf{g}}(\mathbf{y})
\approx
\hat{\mathbf{g}}_s(\mathbf{y})
=
\sum_{n=M+1}^{N}
\mathbf{a}_n(\mathbf{y})
f^n(\mathbf{y}).
\end{equation}

The $M$ algebraic constraints of Eq.~\eqref{eq:CSP_sim} define an implicit approximation of the SIM.~The fast variables $\mathbf{y}_f\in\mathbb{R}^{M}$ can thus be represented as smooth functions of the slow variables $\mathbf{y}_s\in\mathbb{R}^{N-M}$:
\begin{equation}
\label{eq:CSP_slow_fast_relation}
\mathbf{y}_f
=
h(\mathbf{y}_s).
\end{equation}
Near-manifold samples: 
\[
\{
(\mathbf{y}_s^{(k)},\mathbf{y}_f^{(k)})
\}_{k=1}^{T'},
\]
are extracted from the trajectory, by retaining time instances for which the fast modal amplitudes remain below prescribed tolerances.~The manifold map is then identified through regression:
\begin{equation}
\label{eq:manifold_map}
h^\star
\in
\arg\min_{h\in\mathcal{H}}
\frac{1}{T'}
\sum_{k=1}^{T'}
\left\|
\mathbf{y}_f^{(k)}
-
h(\mathbf{y}_s^{(k)})
\right\|^2,
\end{equation}
subject to structural or complexity constraints on the hypothesis space $\mathcal{H}$ \cite{shalev2014understanding}.~Once the SIM has been identified as:
\begin{equation}
\label{eq:SIM_slow_fast_relation}
\mathbf{y}_f = h^\star(\mathbf{y}_s),
\end{equation}
the reduced-order dynamics are obtained by restricting the reconstructed vector field to the manifold:
\begin{equation}
\label{eq:SIM_slow_variables}
\dot{\mathbf{y}}_s
=
\hat{\mathbf{g}}_s
\big(
\mathbf{y}_s,
h^\star(\mathbf{y}_s)
\big).
\end{equation}
To identify an explicit sparse representation of the reduced-order dynamics, manifold-consistent trajectories are constructed by sampling the fast variables from the above regression eq.~\ref{eq:SIM_slow_fast_relation}:
\begin{equation}
\label{eq:SIM_trajectory}
\mathbf{y}^{\mathrm{SIM}}(t_k)
=
\big(
\mathbf{y}_s(t_k),
h^\star(\mathbf{y}_s(t_k))
\big),
\qquad
k=1,\dots,T.
\end{equation}
These trajectories lie exactly on the learned manifold and are used as input for reduced-model identification.

To improve robustness in multi-scale temporal domains, while avoiding explicit numerical differentiation, the reduced-order dynamics are identified using a weak (integral) formulation of SINDy~\cite{messenger2021weak}, in which time derivatives are replaced by integral constraints obtained through smooth test functions and integration by parts. 

Sparse regression with thresholding is then applied to identify a low-complexity reduced model:
\begin{equation}
\label{eq:SIM_reduced_system}
\dot{\mathbf{y}}_s
\approx
\hat{\mathbf{g}}_{\mathrm{red}}(\mathbf{y}_s),
\end{equation}
which is consistent with the learned SIM. Please note, that any other available technique can be used for the reduced-order dynamics identification \cite{schmidt2009distilling, petersen2019deep, udrescu2020ai}.

\subsection{Algorithmic summary}
\label{sec:algorithm}

The complete workflow is summarized in Algorithm~\ref{alg:framework}. The framework consists of:
(1) temporal data augmentation,
(2) sparse vector-field reconstruction,
(3) matrix-free Jacobian estimation,
(4) CSP multiscale diagnostics,
(5) regression-based SIM identification, and
(6) reduced model construction on the learned SIM.

\begin{algorithm}[h!]
\caption{Data-driven SIM identification and reduced-model construction.}
\label{alg:framework}
\begin{minipage}{0.98\linewidth}
\footnotesize
\vspace{1ex}
\begin{mdframed}[style=gray]
\begin{algorithmic}[0]
\State \textbf{Inputs:} time-series data $\{\mathbf{y}(t_k)\}_{k=1}^{T}\subset\mathbb{R}^{N}$, sampling interval $\Delta t$, feature library $\boldsymbol{\Theta}(\mathbf{y})$, sparse-regression settings, Krylov dimension $r\le N$, finite-difference step $\varepsilon$, symbolic regression hypothesis class $\mathcal{H}$, weak-form parameters, regularization constant $\eta>0$.
\State \textbf{Outputs:} smooth surrogate $\hat{\mathbf{g}}(\mathbf{y})$, Jacobian estimates $\{\mathbf{J}_k\}_{k=1}^{T}$, diagnostics $\{\tau_n(t_k), f^n(t_k), Po^n(t_k)\}$, SIM map $h^\star$, reduced model $\hat{\mathbf{g}}_{\mathrm{red}}$.
\end{algorithmic}
\end{mdframed}
\vspace{-1.7mm}
\begin{mdframed}[style=dmd_blue]
\begin{algorithmic}[0]
\State[1] \textbf{Data augmentation via Dynamic Mode Decomposition (DMD).}
Use DMD to reconstruct state trajectories on a refined time grid from sparsely sampled data.
\end{algorithmic}
\end{mdframed}
\vspace{-1.7mm}
\begin{mdframed}[style=sindy_purple]
\begin{algorithmic}[0]
\State[2] \textbf{Smooth vector-field reconstruction.}
Reconstruct a sparse surrogate of the unknown dynamics $\mathbf{g}$ from data using SINDy:
\[
\hat{\mathbf{g}}(\mathbf{y})=\boldsymbol{\Theta}(\mathbf{y})\,\boldsymbol{\Xi},
\]
where $\boldsymbol{\Xi}$ is identified via SR3 or STLSQ.
Select the final surrogate $\hat{\mathbf{g}}(\mathbf{y})$ by threshold sweeps and, if desired, ensemble aggregation.
\end{algorithmic}
\end{mdframed}
\vspace{-1.7mm}
\begin{mdframed}[style=jacobian_orange]
\begin{algorithmic}[0]
\State[3a] \textbf{Define a matrix-free Jacobian--vector product operator.}
For each sample $\mathbf{y}(t_k)$ and direction $\mathbf{v}\in\mathbb{R}^{N}$, approximate
\[
\mathbf{J}(\mathbf{y})\,\mathbf{v}\approx
\big(\hat{\mathbf{g}}(\mathbf{y}+\varepsilon\mathbf{v})-\hat{\mathbf{g}}(\mathbf{y})\big)/\varepsilon.
\]
\State[3b] \textbf{Krylov--Arnoldi Jacobian approximation along the trajectory.}
\For{$k=1,\dots,T$}
\State Initialize a random unit vector $\mathbf{v}_1\in\mathbb{R}^{N}$.
\State \textbf{Arnoldi iteration:} for $j=1,\dots,r$,
\State \hspace{3mm} compute $\mathbf{w}=\mathbf{J}(\mathbf{y}(t_k))\,\mathbf{v}_j$ using Step [3a];
\State \hspace{3mm} orthonormalize $\mathbf{w}$ against $\{\mathbf{v}_1,\dots,\mathbf{v}_j\}$ (Gram--Schmidt);
\State \hspace{3mm} assemble $\mathbf{V}_r\in\mathbb{R}^{N\times r}$ and $\mathbf{H}_r\in\mathbb{R}^{r\times r}$.
\State Construct the Galerkin approximation:
\[
\mathbf{J}_r=\mathbf{V}_r\mathbf{H}_r\mathbf{V}_r^{\mathsf{T}}.
\]
\EndFor
\end{algorithmic}
\end{mdframed}
\vspace{-1.7mm}
\begin{mdframed}[style=csp_teal]
\begin{algorithmic}[0]
\State[4] \textbf{CSP-inspired diagnostics from the estimated Jacobian.}
\For{$k=1,\dots,T$}
\State Compute eigenpairs $\{(\lambda_n,\mathbf{a}_n,\mathbf{b}^n)\}_{n=1}^{N}$ of $\mathbf{J}_k$, developing time scales $\tau_n=|\lambda_n|^{-1}$,
\State modal amplitudes $f^n=\mathbf{b}^n\, \hat{\mathbf{g}}(\mathbf{y})$ and pointer values $Po_i^n=a_n^i\,b_i^n$ for $n=1,\dots,N$.
\State Identify fast/slow variables and infer equilibrium-type indicators.
\EndFor
\end{algorithmic}
\end{mdframed}
\vspace{-1.7mm}
\begin{mdframed}[style=sim_green]
\begin{algorithmic}[0]
\State[5a] \textbf{SIM data extraction and constraint definition.}
Use the diagnostics to define a constraint $\mathbf{r}(\mathbf{y})=\mathbf{0}$.
In QSS-type regimes, represent the SIM as a graph $\mathbf{y}_f=h(\mathbf{y}_s)$ and extract near-manifold samples
$\{(\mathbf{y}_s^{(j)},\mathbf{y}_f^{(j)})\}_{j=1}^{T'}$.
\State[5b] \textbf{Regression-based SIM identification.}
Fit an explicit manifold map $h\in\mathcal{H}$ by regression:
\[
h^\star \in
\arg\min_{h\in\mathcal{H}}
\frac{1}{T'}
\sum_{j=1}^{T'}
\left\|
\mathbf{y}_f^{(j)}
-
h(\mathbf{y}_s^{(j)})
\right\|^2.
\]
In the present work, this regression is performed using PySR, which searches for an explicit analytical expression while balancing prediction accuracy and expression complexity.
\end{algorithmic}
\end{mdframed}
\vspace{-1.7mm}
\begin{mdframed}[style=rom_magenta]
\begin{algorithmic}[0]
\State[6] \textbf{Reduced-order model identification on the learned manifold.}
Construct manifold-consistent trajectories by enforcing $\mathbf{y}_f=h^\star(\mathbf{y}_s)$.
Identify a sparse reduced-order model (restricted to $\mathbf{y}_s$) using weak SINDy:
\[
\dot{\mathbf{y}}_s \approx \hat{\mathbf{g}}_{\mathrm{red}}(\mathbf{y}_s),
\]
and select the final model via sparsity--accuracy trade-off criteria.
\end{algorithmic}
\end{mdframed}
\vspace{1ex}
\end{minipage}
\end{algorithm}

\section{Why regression recovers the slow manifold: an analytically solvable  stochastic system}
\label{section3}

To analytically validate the proposed framework, we consider a two-dimensional stochastic linear system with explicitly enforced time-scale separation. Due to the linearity of the dynamics, the slow invariant manifold (SIM), the fast variable, and the associated reduced model can all be derived analytically, providing a transparent benchmark for validating the CSP diagnostics and the regression-based SIM identification procedure.

Beyond validation, the stochastic system considered admits an exact \emph{Mori--Zwanzig} (MZ) interpretation that clarifies \emph{why} the regression-based reduction succeeds \cite{mori1965transport,zwanzig1961memory,chorin2000optimal}. Treating the slow coordinate as the resolved variable and the fast coordinate as the unresolved variable, we show that: (i) the regression-based SIM is precisely the \emph{Mori projection} of the fast variable onto the slow one under the stationary measure, i.e. the conditional expectation $\mathbb{E}[y\mid x]$; (ii) the vector field of the reduced model is the \emph{Markovian streaming term} of the corresponding generalized Langevin equation (GLE); and (iii) the CSP statement that the fast modal amplitude is exhausted ($f^1\approx 0$) is the same dynamical fact as the \emph{collapse of the MZ memory kernel} to a Markovian limit. The CSP slow subspace, the geometric SIM, and the statistical Mori projection thus coincide at leading order; i.e. three readings of one object. This also supplies a statistical-mechanical justification, complementary to the CSP one, for the framework's empirical robustness to noise exploited in the Applications section. A fully self-contained, pedagogical derivation of all MZ statements made below, which requires no prior familiarity with projection-operator methods, is provided in S1 Appendix, Section~\ref{app:MZ}.

\paragraph{System definition and time-scale separation}

Let the state vector be $\mathbf{z}_t=(x_t,y_t)^\top$.~The dynamics are governed by the stochastic differential equation:
\begin{equation}
d\mathbf{z}_t
=
\mathbf{J}\mathbf{z}_t\,dt
+
d\mathbf{W}_t,
\end{equation}
where $\mathbf{W}_t$ denotes a two-dimensional Wiener process and:
\begin{equation}
\mathbf{J}
=
\begin{pmatrix}
-a & b\\
c & -d
\end{pmatrix}.
\end{equation}
where $a,b,c,d$ are positive constant parameters.~To enforce a singularly perturbed regime in which $y$ evolves on the fast time scale and $x$ evolves on the slow time scale, we impose:
\begin{equation}
d \gg a,b,c.
\end{equation}

The eigenvalues of $\mathbf{J}$ are:
\begin{equation}
\lambda_{1,2}=\frac{-(a+d)\pm\sqrt{(a-d)^2+4bc}}{2}.
\end{equation}
and, under the assumption $d \gg a, b, c$, the eigenvalues simplify to:
\begin{equation}
|\lambda_1| \approx |-d| = d \gg |\lambda_2|=|-a + \frac{bc}{d}| \approx a. 
\end{equation}
This explicitly proves a massive separation of dissipative time scales ($\tau_1 \ll \tau_2$), isolating $\tau_1$ as the highly dissipative time scale.~Employing the quasi steady-state approximation, $\dot{y} = cx - dy \approx 0$, the SIM reads $y=
(c/d)x$ and the corresponding slow model is $\dot{x}=\lambda_2 x$.~The remainder of this section recovers these objects in three independent ways and identifies them with the MZ decomposition.

\paragraph{Fast/slow splitting and the generalized Langevin equation}

We designate $x$ as the \emph{slow} (resolved) variable and $y$ as the \emph{fast} (unresolved) variable, and equip the space of observables with the stationary inner product $\langle A,B\rangle=\mathbb{E}_{\mathrm{ss}}[AB]$ induced by the invariant Gaussian measure. Let $\mathcal{P}$ denote orthogonal projection onto the resolved subspace $\mathrm{span}\{x\}=x$ and $\mathcal{Q}=\mathcal{I}-\mathcal{P}$ its complement.~Eliminating $y$ exactly by variation of parameters (see S1 Appendix, Section~\ref{app:MZ}) yields a closed, non-Markovian equation for the resolved variable of Mori form:
\begin{equation}
\dot{x}(t)=-a\,x(t)+\int_0^t K(t-s)\,x(s)\,ds+F(t),
\label{eq:GLE_main}
\end{equation}
with exact memory kernel and fluctuating (orthogonal-dynamics) force:
\begin{equation}
K(t)=bc\,e^{-dt},
\qquad
F(t)=b\,e^{-dt}y(0)+b\!\int_0^t e^{-d(t-s)}\eta_y(s)\,ds+\eta_x(t),
\label{eq:kernel_force_main}
\end{equation}
where $\eta_x,\eta_y$ are the components of the white noise. The streaming frequency is fixed by the stationary statistics:
\begin{equation}
\Omega=\frac{\langle\dot{x},x\rangle}{\langle x,x\rangle}=-a+b\,\frac{\Sigma_{xy}}{\Sigma_{xx}}=-a+\frac{bc}{d}=\lambda_2.
\label{eq:Omega_main}
\end{equation}
Both the memory kernel and the fluctuating force relax on the \emph{fast} time scale $\tau_{\mathrm{mem}}=d^{-1}=\tau_1$.~In the singular limit $d\gg a,b,c$ the kernel acts Markovianly, $\int_0^t e^{-d(t-s)}x(s)\,ds\to x(t)/d$, so that after the initial fast layer Eq.~\eqref{eq:GLE_main} collapses to:
\begin{equation}
\dot{x}=\Big(-a+\frac{bc}{d}\Big)x=\lambda_2\,x.
\label{eq:reduced_collapse_main}
\end{equation}
Thus the reduced slow dynamics is the Markovian streaming term of the GLE, and the exhaustion of the fast CSP mode ($f^1\approx 0$; see S3 Appendix, Section~\ref{app:csp_2d}) is identical to the collapse of the MZ memory and noise terms. Equivalently, $\Omega=-a+\int_0^\infty K(t)\,dt=-a+bc/d$: the effective slow rate is the bare rate dressed by the integrated memory.

\paragraph{Identification of the fast variable}

The fast variable will be identified using the CSP pointer, as follows; detailed steps are presented in S1 Appendix, Section~\ref{app:fast_variable}.~The right eigenvector associated with the fast eigenvalue $\lambda_1$ of $\mathbf{J}$ satisfies:
\[
\mathbf{J}\mathbf{a}_1
=
\lambda_1\mathbf{a}_1.
\]
Under the leading-order approximation $d\gg a,b,c$, the corresponding right eigenvector is:
\begin{equation}
\mathbf{a}_1
\approx
\begin{bmatrix}
b\\
-d
\end{bmatrix},
\end{equation}
while the associated left eigenvector, appropriately normalized to satisfy biorthonormality ($\mathbf{b}^1\mathbf{a}_1 = 1$), is:
\begin{equation}
\mathbf{b}^{1}
\approx
\frac{1}{d^2}
\begin{bmatrix}
c & -d
\end{bmatrix}.
\end{equation}

The CSP pointer associated with the fast mode is then:
\begin{equation}
Po^1
=
\operatorname{diag}
(\mathbf{a}_1\mathbf{b}^{1})
=
\begin{bmatrix}
bc/d^2\\
1
\end{bmatrix}
\approx
\begin{bmatrix}
0\\
1
\end{bmatrix}.
\end{equation}

The pointer therefore localizes the fast dissipative mode almost entirely on the second variable, correctly identifying $y$ as the fast variable. In MZ language, the deterministic CSP fast projector $\mathbf{a}_1\mathbf{b}^{1}$ coincides at leading order with the statistical projector $\mathcal{Q}$: both annihilate the slaved part of $y$ and isolate the fast fluctuation that carries the memory and fluctuating force of Eq.~\eqref{eq:kernel_force_main}.

\paragraph{The slow invariant manifold as a Mori projection}
The expectation of the state vector $\boldsymbol{\mu}(t) = \mathbb{E}[\mathbf{z}_t]$ evolves according to the deterministic system $\dot{\boldsymbol{\mu}} = \mathbf{J} \boldsymbol{\mu}$. Because $\mathbf{J} $ is stable ($\text{tr}(\mathbf{J} ) < 0$ and $\text{det}(\mathbf{J} ) > 0$), the means asymptotically decay to zero:
\begin{equation}
\lim_{t \to \infty} \mu_x(t) = 0, \quad \lim_{t \to \infty} \mu_y(t) = 0
\end{equation}

The stationary covariance matrix:
\[
\mathbf{\Sigma}
=
\mathbb{E}
\left[
(\mathbf{z}_t-\boldsymbol{\mu}_t)
(\mathbf{z}_t-\boldsymbol{\mu}_t)^\top
\right]
\]
satisfies the continuous-time Lyapunov equation:
\[
\mathbf{J}\mathbf{\Sigma}
+
\mathbf{\Sigma}\mathbf{J}^\top
+
\mathbf{I}
=
0.
\]

Using dominant-balance arguments under the assumption $d\gg a,b,c$ (see S1 Appendix, Section~\ref{app:covariance_dynamics} for details), the stationary covariance components satisfy:
\begin{equation}
\Sigma_{xx}
\approx
\frac{1}{2a},
\qquad
\Sigma_{yy}
\approx
\frac{1}{2d},
\qquad
\Sigma_{xy}
\approx
\frac{c}{2ad}.
\end{equation}

For the zero-mean jointly Gaussian stationary state, the Mori projection of the unresolved variable onto the resolved one is exactly the conditional expectation:
\begin{equation}
\boxed{\;\mathcal{P}y=\frac{\langle y,x\rangle}{\langle x,x\rangle}\,x=\frac{\Sigma_{xy}}{\Sigma_{xx}}\,x=\frac{c}{d}\,x=\mathbb{E}[y\mid x].\;}
\label{eq:Py_main}
\end{equation}

After identifying $y$ as the fast variable, the SIM is estimated through linear regression of the form:
\[
y
=
\beta_1 x+\beta_0+\epsilon.
\]

The regression coefficients are strictly defined by the system's theoretical moments. The intercept is $\beta_0 = \mu_y - \beta_1 \mu_x$.
Since $\mu_x = \mu_y = 0$ at steady state, the statistical intercept vanishes: $\beta_0 = 0$.

Using the standard regression estimator \cite{fahrmeir2013regression}, $\beta_1=\mathrm{Cov}(x,y)/\mathrm{Var}(x)=\Sigma_{xy}/\Sigma_{xx}$ and substituting the asymptotic covariance approximations yields:
$\beta_1 \approx (c/2ad)/(1/2a)=c/d$.~Additional details regarding the associated Jacobian inverse problem and the asymptotic recovery of the dominant dynamical parameters from the stationary covariance matrix are provided in S1 Appendix, Section~\ref{app:inverse_problem}.~Comparing with Eq.~\eqref{eq:Py_main}, the regression slope is precisely the Mori projection coefficient, and the learned manifold $y=(c/d)x$ is the Mori projection $\mathcal{P}y$.~The data-driven regression slope is mathematically identical to the coefficient of the deterministic SIM, confirming that the noise cloud natively aligns its principal axis along the slow dynamical subspace.

The fluctuation discarded by the regression is exactly the orthogonal component $\mathcal{Q}y=y-(c/d)x$, whose projection onto $x$ vanishes by construction and whose residual variance,
\begin{equation}
\langle\mathcal{Q}y,\mathcal{Q}y\rangle
=\Sigma_{yy}-\frac{\Sigma_{xy}^2}{\Sigma_{xx}}
=\frac{1}{2d}\Big(1-\frac{c^2}{ad}\Big)\approx\frac{1}{2d},
\label{eq:Qy_main}
\end{equation}
is the $\mathcal{O}(d^{-1})$ fast noise that lives entirely off the manifold.~This is the rigorous content of the alignment statement above: the regression returns the conditional mean (the manifold), which is the minimum-variance estimator and is therefore natively insensitive to the orthogonal fast fluctuations.

In the deterministic CSP framework, setting the fast dynamics to be exhausted ($\mathbf{b}^{1}\dot{\mathbf{z}_t} =\dot{y} = cx - dy \approx 0$) yields the exact geometric equation for the SIM: $y = (c/d)x$, and the slow right eigenvector is $\mathbf{a}_2\propto(1,\,c/d)^\top$.~The CSP slow eigenspace, the regression manifold, and the Mori projection coincide at leading order in $d^{-1}$.



\section{Applications}

Previous section establishes the link between regression and the SIM in the
one setting where it can be computed exactly: a linear system driven by white noise.
Here, least-squares regression of $y_f$ on $y_s$ under the stationary measure is the
Mori projection onto linear functions of $y_s$. For the two-dimensional example (slow
$x$, fast $y$), let $\delta = a/d\ll1$ denote the ratio of fast to slow timescales at
leading order, and normalize the slow right eigenvector as
$\mathbf{a}_2 = (1, m)^\top$, so that the SIM is $y = m\, x$. For $c \neq 0$, the
regression slope $\beta_1 = \Sigma_{xy}/\Sigma_{xx}$ then satisfies
\begin{equation}
   \frac{\text{Regression slope}}{\text{SIM slode}} =\frac{\beta_1}{m} = 1 + \delta\left(\frac{b}{c} - 2\right) + \mathcal{O}(\varepsilon^2)
  \label{eq:beta_over_sim}
\end{equation}
(S1 Appendix, Section~\ref{SI:convergence}). Regression therefore recovers the SIM at leading order, with
a relative discrepancy of order $\delta$. This correction, however, is not a property
of the manifold: the term $b/c$ is noise on the fast variable transmitted to $x$ through
the back-coupling (S1 Appendix), and in the limiting case $c = 0$ the SIM is $y = 0$, yet
$\beta_1 \neq 0$. For nonlinear systems the corresponding object is Zwanzig's projection
$\mathbb{E}[y_f \mid y_s]$, which the framework estimates by regression over a nonlinear
hypothesis space $\mathcal{H}$. Exhaustion of the fast modal amplitudes and reduction of
the memory kernel to its Markovian limit are both consequences of the spectral gap;.

The benchmarks below are deterministic, and here the argument is more direct. Once the
fast amplitudes fall below the exhaustion tolerance, trajectories lie within a distance
of the attracting slow manifold that is determined by that tolerance and decays exponentially
in time. The regression therefore targets the SIM itself rather than a leading-order
approximation of it, and no noise-induced correction of the kind in
Eq.~\eqref{eq:beta_over_sim} arises. Its accuracy is limited by the exhaustion tolerance,
the capacity of $\mathcal{H}$, and the part of the manifold the data cover: the SIM is
identified only where trajectories pass. The fast--slow partition and the exhaustion
windows, by contrast, come from the Jacobian of the reconstructed vector field, so their
reliability depends on the surrogate rather than on the argument above. Because data
lying on the SIM carry information about the fast rates only through small deviations
from it, the ordering of timescales and the CSP pointers are more robust than the
timescale values.  
The four following benchmarks probe different parts of this argument: a linear SIM (2D model), a
switch between regimes with an exchange of fast variable (Michaelis - Menten), a
near-linear SIM with a higher-order correction (biofilm), and rational SIMs in two
successive regimes (TMDD). For each  benchmark system considered, the framework is used to identify fast and slow dynamics through CSP diagnostics, construct data-driven approximations of the SIM, and identify reduced models governing the slow evolution on the learned SIM. In all cases considered, a dense data set is employed, so that the issue of sparsity is not addressed; this issue is discussed in Section \ref{sec:five}.~Numerical settings, parameter values, initial conditions, and sampling details are reported in S1 Appendix, Section~\ref{app:numerical_settings}. Additional CSP diagnostic results, including reference, downsampled, and noisy-data cases, are reported in S2 Appendix, Section~\ref{app:csp_diagnostics}.

The analysis of each benchmark system follows a uniform three-stage structure.~First, we establish the system kinetics and investigate its local time-scale spectrum and Jacobian eigenstructure using matrix-free Krylov--Arnoldi diagnostics to track the fast-slow variable partition.~Second, we extract near-manifold trajectory snapshots, in which the fast dissipative modal amplitudes are exhausted ($f^m \approx 0$), to map the geometry of the SIM via symbolic or algebraic regression.~Finally, we restrict the state space to this learned manifold and employ Weak SINDy to discover a low-complexity, explicit reduced-order model.~The accuracy of each stage is validated by directly comparing our data-driven manifold maps and reconstructed vector fields against their analytically derived singular perturbation counterparts.

\subsection{Illustrative 2D linear model}

\subsubsection{Model and multi-scale characteristics}

\begin{figure}[h!]
\centering
\includegraphics[width=0.75\linewidth]{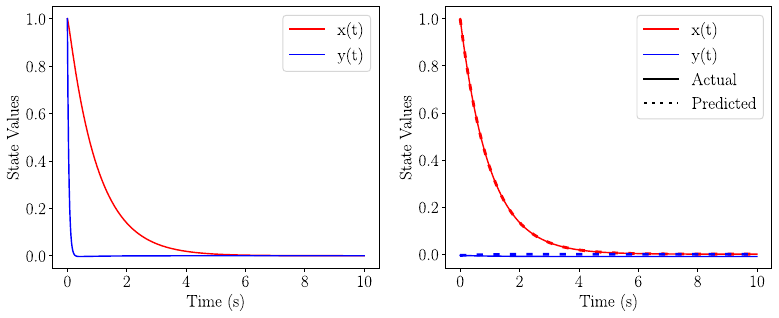}
\caption{Reference trajectories and data-driven model performance for the 2D linear system.~Left: Trajectories from the full system.~Right: comparison of trajectories on the SIM obtained by the analytically derived reduced model (actual) and the data-driven model (predicted).}
\label{fig:linear2D_unified}
\end{figure}

We first consider a two-dimensional linear system, of the form used in Section~\ref{section3}, serving as an analytically tractable benchmark:
\begin{equation}
\begin{aligned}
\dot x &= -x + 0.5\,y,\\
\dot y &= -0.1\,x - 20\,y.
\end{aligned}
\label{eq:linear-system}
\end{equation}
The constant Jacobian $\mathbf{J}$ for this system yields time-independent eigenvalues that reveal a pronounced, stationary time scale separation.~The local CSP pointer diagnostic localizes the fast dissipative mode entirely on the coordinate $y$, while the slow dynamics are parameterized by $x$, driving the trajectories into rapid alignment with the slow linear subspace ($y\approx 0$, see left panel of Fig~\ref{fig:linear2D_unified}).~Full reference diagnostics are provided in S3 Appendix, Section ~\ref{app:csp_2d_reference_diagnostics}.

\subsubsection{SIM identification}

Table \ref{tab:error_results_2D} displays the mean square error (MSE) and root mean square error (RMSE) of the fast variable $y$, as computed  by the analytically derived reduced model and the data-driven model.~The results demonstrate an excellent agreement of the two models.

\begin{table}[h!]
\centering
\caption{Error for the fast variable in the 2D linear model.}
\label{tab:error_results_2D}
\begin{tabular}{c cc}
\toprule
\textbf{Fast variable} & \textbf{MSE} & \textbf{RMSE} \\
\midrule
$y$ & $3.38 \times 10^{-9}$ & $5.82 \times 10^{-5}$ \\
\bottomrule
\end{tabular}
\end{table}

\subsubsection{Reduced-order dynamics}

Table~\ref{tab:comparison_2D} displays the analytic form of the reduced model derived analytically and of the data-driven model derived with the proposed framework.~The two models are in excellent agreement.

\begin{table}[h!]
\centering
\caption{Analytical and identified models for the 2D linear system.}
\label{tab:comparison_2D}
\begin{tabular}{c c}
\toprule
\textbf{Analytical reduced model} & \textbf{Identified model} \\
\midrule
$\dot{x}=-1.0025x$ & $\dot{x}=-1.0030x$ \\
\bottomrule
\end{tabular}
\end{table}

\FloatBarrier
\subsection{Michaelis--Menten model}

\subsubsection{Model and multiscale characteristics}

\begin{figure}[h!]
\centering
\includegraphics[width=1.\linewidth]{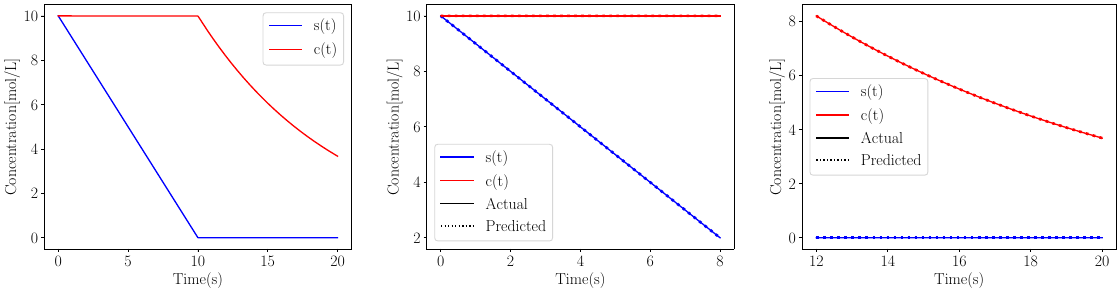}
\caption{Reference trajectories and data-driven model performance for the Michaelis--Menten model. Left:~Trajectories of $c(t)$ and $s(t)$ from the full model.~Middle: Comparison of analytical and data-driven trajectories within the early-time sQSSA regime.~Right: Comparison of analytical (actual) and data-driven (predicted) trajectories within the late-time rQSSA regime.}
\label{fig:mm_unified}
\end{figure}

We next apply the framework to the Michaelis--Menten model \cite{PATSATZIS2019108220}:
\begin{equation}
\begin{aligned}
    \dot{s} &= -k_{1f}(e_{0}-c)s + k_{1b}c, \\
    \dot{c} &= k_{1f}(e_{0}-c)s - (k_{1b}+k_{2})c,
\end{aligned}
\label{eq:MM}
\end{equation}
where $s(t)$ denotes the substrate concentration, $c(t)$ denotes the enzyme-substrate complex concentration, and $e_0$ is the total enzyme concentration. The full-order unreduced reference trajectories are illustrated in panel A of Fig~\ref{fig:mm_unified}.

The time scale spectrum exhibits a clear separation between one fast and one slow dissipative mode over practically the full phase and parameter spaces \cite{PATSATZIS2019108220}.~For the specific parameters used here, CSP diagnostics along the SIM reveal that the fast variable changes identity across the dynamics: in the early regime the standard quasi-steady-state approximation (sQSSA) applies, according to which \(c\) is fast, whereas in the late regime  the reverse quasi-steady-state approximation (rQSSA) applies, according to which \(s\) is fast \cite{PATSATZIS2023}.~The corresponding developing time scales, modal amplitudes, and pointer diagnostics are reported in S3 Appendix, Section~\ref{app:csp_MM}.

\subsubsection{SIM identification}

Because the Michaelis--Menten system exhibits distinct sQSSA and rQSSA regimes along the trajectory on the SIM, the portion of the SIM that refers to each of the two QSSAs is identified through data-driven regression on the corresponding partitioned dataset. 

Using the CSP pointer, it is recognized that the sQSSA  is valid when \(t < 10\) and the rQSSA is valid when \(t > 10\).~In the sQSSA regime the manifold is learned by regressing \(c\) against \(s\), resulting in the graph $c=h_{\mathrm{sQSSA}}(s)$, while in the rQSSA regime the manifold is obtained by regressing \(s\) against \(c\) to yield $s=h_{\mathrm{rQSSA}}(c)$.~This partition enables direct comparison against the classical analytical forms of the SIM:
\begin{equation}
c = \frac{k_{1f}e_0s}{k_{1f}s+k_{1b}+k_2}, \qquad s = \frac{k_{1b}c}{k_{1f}(e_0-c)}. 
\end{equation}
The errors of the fast variable for both learned manifolds are reported in Table~\ref{tab:error_results_MM}.~Once more, excellent agreement is demonstrated.

\begin{table}[h!]
\centering
\caption{Errors for the fast variables in the Michaelis--Menten model.}
\label{tab:error_results_MM}
\begin{tabular}{c cc}
\toprule
\textbf{Fast variable} & \textbf{MSE} & \textbf{RMSE} \\
\midrule
$c \text{ (sQSSA)}$ & $3.13 \times 10^{-11}$ & $5.59 \times 10^{-6}$ \\
$s \text{ (rQSSA)}$ & $9.56 \times 10^{-13}$ & $9.78 \times 10^{-7}$ \\
\bottomrule
\end{tabular}
\end{table}

\subsubsection{Reduced-order dynamics}

Separate data-driven models were generated with the proposed framework on the basis of the  sQSSA and rQSSA datasets.~These models are compared with the analytically derived reduced models in Table~\ref{tab:comparison_MM_models}, demonstrating an excellent agreement.~This agreement is reflected in identical profiles, shown in the two right panels of Fig~\ref{fig:mm_unified}.

\begin{table}[h!]
\centering
\caption{Analytical and identified models for the Michaelis--Menten system.}
\label{tab:comparison_MM_models}
\begin{tabular}{c c c}
\toprule
\textbf{Regime} & \textbf{Analytical reduced model} & \textbf{Identified model} \\
\midrule
sQSSA & $\dot{s} = -0.100\,c$ & $\dot{s} = -0.100\,c$ \\
rQSSA & $\dot{c} = -0.100\,c$ & $\dot{c} = -0.100\,c$ \\
\bottomrule
\end{tabular}
\end{table}

\FloatBarrier
\subsection{Biofilm model}

\subsubsection{Model and multi-scale characteristics}

We consider the three-dimensional nonlinear biofilm growth model \cite{gardedifferential}:
\begin{equation}
\begin{aligned}
\dot G_p &= k_1\,G_E\,G_p - k_4\,G_p - k_2\,A\,G_p,\\
\dot A   &= -k_3\,A + k_5\,G_i,\\
\dot G_i &= k_4\,G_p - k_5\,G_i.
\end{aligned}
\label{eq:biofilm}
\end{equation}

\begin{figure}[h!]
\centering
\includegraphics[width=0.75\linewidth]{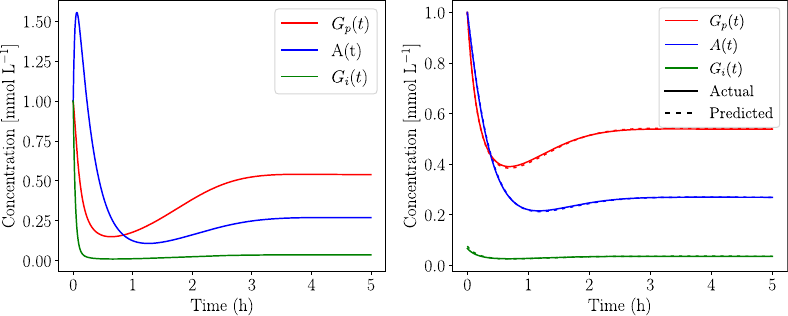}
\caption{Reference trajectories and data-driven model performance for the biofilm model. Left: Trajectories from the full system.~Right: comparison of trajectories on the SIM obtained by the analytically derived reduced model (actual) and the data-driven model (predicted).}
\label{fig:biofilm_time_unified}
\end{figure}

The dynamics of this system exhibit a persistent spectral gap, separating one fast dissipative mode from two slow modes.~The CSP pointer consistently associates $G_i$ with the fast exhausted mode, confirming that $(G_p, A)$ act as the  slow coordinates parameterizing the low-dimensional flow.~The trajectories displayed in the left panel of Fig~\ref{fig:biofilm_time_unified}, demonstrate the multi-scale behavior, by the steep initial profiles and the smooth later ones.~Full CSP diagnostics  are detailed in S3 Appendix, Section~\ref{app:csp_biofilm_reference_diagnostics}.

\subsubsection{SIM identification}

\begin{figure}[h!]
\centering
\includegraphics[width=0.6\linewidth]{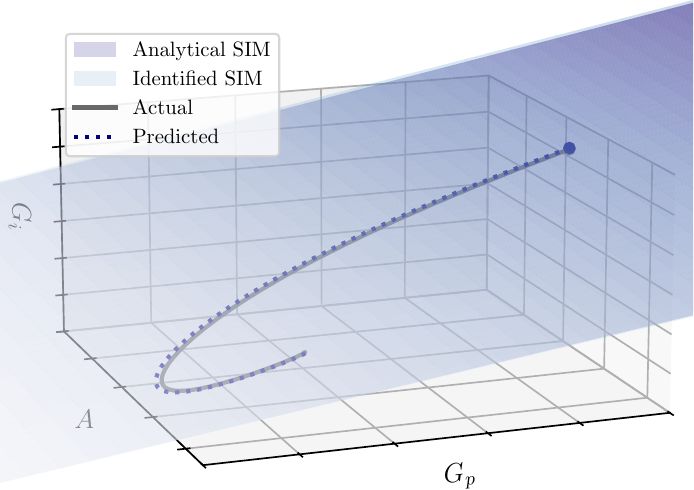}
\caption{SIMs and trajectories comparison for the biofilm model.~The data-driven identified SIM is compared to the analytical SIM ($G_i$-based QSSA).~Trajectories from the data-driven model (predicted) and from the analytical reduced model (actual), initiated on the related SIMs, are also displayed.}
\label{fig:biofilm_geom_unified}
\end{figure}

Symbolic regression evaluated over the near-manifold snapshot data produces the following explicit algebraic surface mapping:
\begin{equation}
\label{eq:biofilm_SIM}
G_i \approx G_p\Big(\alpha+\beta\,A\Big)+ \gamma,
\end{equation}
where $(\alpha,\beta,\gamma)$ are the discovered numerical coefficients; see S3 Appendix, Section~\ref{app:csp_biofilm_reference_diagnostics}.~Figure~\ref{fig:biofilm_geom_unified} compares the learned SIM  with the analytical SIM:
\begin{equation}
\label{eq:biofilm_SIMx}
G_i \approx G_P \frac{k_4}{k_5} \left[ 1 - \frac{k_1}{k_5} G_E + \frac{k_4}{k_5} + \frac{k_2}{k_5} A \right],
\end{equation}
which was derived by considering higher order accuracy of the QSSA of $G_i$ \cite{goussis2012quasi}.~An excellent agreement is demonstrated between the two SIMs; see S3 Appendix, Section~\ref{app:csp_biofilm_reference_diagnostics} for details.~The errors on the fast variable $G_i$ reported in Table~\ref{tab:error_results_Biofilm} provide a quantitative assessment of this agreement.

\begin{table}[h!]
\centering
\caption{Errors for the fast variable in the biofilm model.}
\label{tab:error_results_Biofilm}
\begin{tabular}{c cc}
\toprule
\textbf{Fast variable} & \textbf{MSE} & \textbf{RMSE} \\
\midrule
$G_i$ & $1.48 \times 10^{-6}$ & $1.22 \times 10^{-3}$ \\
\bottomrule
\end{tabular}
\end{table}
\subsubsection{Reduced-order dynamics}

The data-driven model generated by the proposed framework, was shown to be in an excellent agreement with the one derived analytically, on the basis the leading order accurate QSSA of $G_i$.~This is shown in Table \ref{tab:comparison_biofilm_models}, where analytic expressions of the two models are displayed.~This agreement is also demonstrated in Fig.~\ref{fig:biofilm_geom_unified}, where trajectories from the analytical reduced and the data-driven models are displayed.~These trajectories initiate from the corresponding SIMs and remain there.

\begin{table}[h!]
\centering
\caption{Analytical and identified models for the biofilm system.}
\label{tab:comparison_biofilm_models}
\begin{tabular}{c c}
\toprule
\textbf{Analytical reduced model} & \textbf{Identified model} \\
\midrule
$\dot{G}_p=1.426G_p-5.300G_pA$ & $\dot{G}_p=1.426G_p-5.300G_pA$ \\
$\dot{A}=2.000G_p-4.000A$ & $\dot{A}=1.937G_p-3.892A$ \\
\bottomrule
\end{tabular}
\end{table}

\FloatBarrier
\subsection{Target-mediated drug disposition model}

\subsubsection{Model and multiscale characteristics}

\begin{figure}[h!]
\centering
\includegraphics[width=\linewidth]{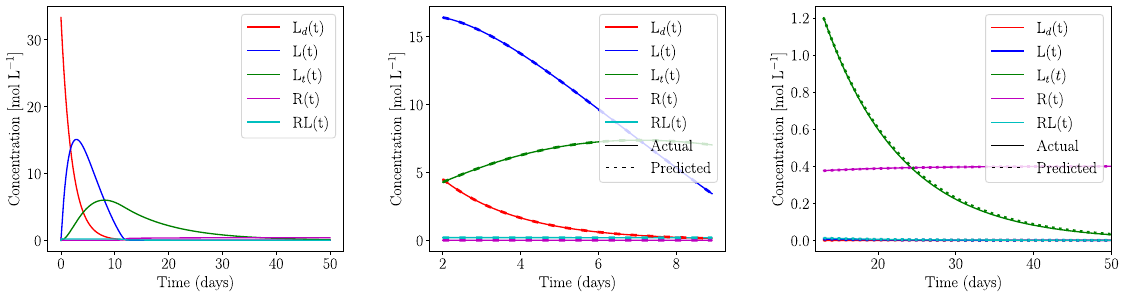}
\caption{Reference trajectories and data-driven model performance for the TMDD model. Left: Trajectories from the full system.~Right: comparison of trajectories on the SIM obtained by the analytically derived reduced model (actual) and the data-driven model (predicted).}
\label{fig:tmdd_3x1_unified}
\end{figure}

Finally, we apply the framework to a complex, five-dimensional target-mediated drug disposition (TMDD) model characterized by highly nonlinear binding kinetics \cite{michalaki2018asymptotic}:
\begin{equation}
\begin{aligned}
&\dot L_d = -k_a L_d,\\
&\dot L_{~}    = k_a L_d - k_{el} L - k_{pt} L  - k_{on} L\,R + k_{off} RL + k_{tp} L_t,\\
&\dot L_t = k_{pt} L - k_{tp} L_t,\\
&\dot R_{~}   \ = -k_{on} L\,R + k_{off} RL + k_{syn} - k_{deg} R,\\
&\dot{RL} = k_{on} L\,R - k_{off} RL - k_{int} RL.
\end{aligned}
\label{eq:TMDD}
\end{equation}
Here, $L_d$ represents the ligand in the donor compartment, $L$ is the free ligand concentration, $L_t$ is the transit ligand, $R$ is the receptor concentration, and $RL$ is the ligand-receptor complex.~Indicative trajectories obtained from the full model are displayed in the left panel of Fig~\ref{fig:tmdd_3x1_unified}.

The dynamics of the TMDD system exhibits an evolving fast/slow time scale spectrum.~During the initial window ($0 \le t < 10$), a time scale gap between two fast dissipative and three slow time scales emerges, in which the CSP pointer identifies the receptor $R$ and the complex $RL$ as the fast variables.~This finding frames $(L_d, L, L_t)$ as the three slow governing variables during this epoch.~At later times ($10<t$), the fast-slow partition shifts to three fast dissipative time scales and two slow ones.~Now the set of slow variables is $(L_d, L_t)$.~Details on the diagnostic outputs are reported in S3 Appendix, Section~\ref{app:csp_TMDD_reference_diagnostics}.

\subsubsection{SIM identification}

Given the two different fast/slow time scale separations, the analysis was carried-out separately in the two regimes.~Considering the initial multi-scale regime, near-manifold data are harvested from the window $2 \le t \le 9$, where the two fast modes  are considered exhausted.~The SIM  is parameterized as a graph mapping the two fast variables to the three slow ones:
\[
R=h_R(L_d,L,L_t),
\qquad
RL=h_{RL}(L_d,L,L_t).
\]
Similarly, for the second multi-scale regime data are harvested from the domain $13 \le t$, where now the three fast modes are considered exhausted.~The SIM  is now parameterized as a graph mapping the three fast variables to the ttwo slow ones:
\[
R=h_R(L_d,L_t),
\qquad
RL=h_{RL}(L_d,L_t).
\qquad
L=h_{RL}(L_d,L_t).
\]
The accuracy of the two data-driven learned SIMs is assessed by comparing trajectories generated by the analytical reduced model that was constructed on the basis of the QSSA for the fast variables with the trajectories generated by the data-driven learned model, as shown in  Table~\ref{tab:error_results_TMDD}.~The errors reported in Table~\ref{tab:error_results_TMDD} demonstrate an excellent agreement of the two SIMs.
\begin{table}[h!]
\centering
\caption{Prediction errors for the fast variables in the TMDD model.}
\label{tab:error_results_TMDD}
\begin{tabular}{cc cc}
\toprule
 
& \textbf{Fast variable} & \textbf{MSE} & \textbf{RMSE} \\
\midrule
\multirow{2}{*}{Regime 1}&	$R$  	& $3.8 \times 10^{-14}$ & $1.95 \times 10^{-7}$ \\
&	$RL$	& $1.14\times 10^{-10}$ & $1.07 \times 10^{-5}$ \\
\midrule
\multirow{3}{*}{Regime 2}&	$L$ 		& $6.99 \times 10^{-10}$ & $2.64 \times 10^{-5}$ \\
&	$R$  	& $5.18 \times 10^{-8}$ & $2.28 \times 10^{-4}$ \\
&	$RL$ 	& $4.66 \times 10^{-9}$ & $6.83 \times 10^{-5}$ \\
\bottomrule
\end{tabular}
\end{table}

\subsubsection{Reduced-order dynamics}

Table~\ref{tab:comparison_TMDD_models1} displays expressions for the analytical QSSA-based reduced model and for the data-driven learned model, derived on the basis of the available data during the period $2 \le t \le 9$.~The findings displayed in Table~\ref{tab:comparison_TMDD_models1} demonstrate a very good agreement.~This is reflected in the trajectories shown in the middle panel of Fig.~\ref{fig:tmdd_3x1_unified} and in the errors displayed in Table ~\ref{tab:error_results_TMDD}.

\begin{table}[h!]
\centering
\caption{Analytical and identified TMDD  models for the first regime.}
\label{tab:comparison_TMDD_models1}
\begin{tabular}{c | c}
\toprule
\textbf{Analytical reduced model} & \textbf{Identified model} \\
\midrule
$\dot{L}_d=-0.500L_d$ 
&
$\dot{L}_d=-0.500L_d$ 
\\
$\dot{L}=0.500L_d-0.150L+0.100L_t-\tfrac{20L}{1+10L}$ 
&
$\dot{L}=0.499L_d-0.150L+0.098L_t-0.999\tfrac{20L}{1+10L}$ 
\\
$\dot{L}_t=0.100L-0.100L_t$ 
&
$\dot{L}_t=0.100L-0.100L_t$ 
\\
\bottomrule
\end{tabular}
\end{table}

Table~\ref{tab:comparison_TMDD_models2} displays the expression for the analytical QSSA-based reduced model.~The data-driven learned model, obtained on the basis of the available data during the period $13 \le t$, is not stated, since it is very complex.~The reason why it deviates from the one shown in Table~\ref{tab:comparison_TMDD_models2} is because the library of Weak SINDy that  is employed here does not include the functional forms in the analytic model.~Nevertheless, the accuracy provided by the data-drivem model is excellent, as shown in the right panel of Fig.~\ref{fig:tmdd_3x1_unified} and in the errors displayed in Table ~\ref{tab:error_results_TMDD}.

\begin{table}[h!]
\centering
\caption{Analytical TMDD  model for the second regime}
\label{tab:comparison_TMDD_models2}
\begin{tabular}{c}
\toprule
\textbf{Analytical reduced model} \\
\midrule
$\dot{L}_d = -0.500 L_d$ \\[2ex]
$\dot{L}_t = \dfrac{3.97 + L_d + 0.2 L_t \pm \sqrt{(3.97 + L_d + 0.2 L_t)^2 + 0.12 L_d + 0.024 L_t}}{6} - 0.100 L_t$\\
\bottomrule
\end{tabular}
\end{table}

\section{Impact of data sparsity and noise on SIM estimation}
\label{sec:five}
\subsection{Data augmentation using Dynamic Mode Decomposition}
\label{sec:DMD}

In many practical applications, time-series data are sparse or irregularly sampled, limiting accurate system identification, Jacobian estimation, and timescale analysis.~To address this, DMD with Hankel pre-processing is used to augment the temporal resolution of the data prior to applying the proposed framework \cite{brunton2017chaos}.~The original trajectories generated in the context of this study are first downsampled to mimic sparse observations.~A Hankel time-delay embedding is then constructed to enrich the temporal information and improve the representation of nonlinear dynamics.~DMD is subsequently applied to reconstruct dense and dynamically consistent trajectories on a finer temporal grid.~The reconstructed data are then processed by the proposed framework. 

Rather than treating DMD purely as a numerical interpolation tool to fill temporal gaps, this step can be interpreted through an operator-theoretic lens aligned with Koopman theory  \cite{brunton2017chaos}.~The DMD-based reconstruction acts as a spectral filter that projects the sparse, unconditioned observations onto a low-dimensional linear state-space sequence.~This pre-conditioning step suppresses high-frequency stochastic variations while preserving the dominant invariant subspace of the underlying system, explaining why subsequent sparse vector-field identification via SINDy yields clean, interpretable terms even when initialized from severely degraded datasets.


\subsubsection{Results}
\begin{figure}[h!]
\centering
\includegraphics[width=\linewidth]{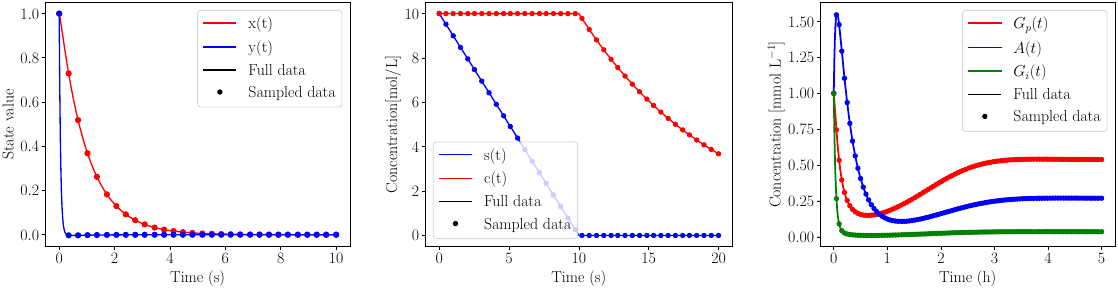}
\caption{
Full trajectories together with the downsampled observations used for DMD-based data augmentation.
Left: 2D linear system with 30 observations.
Middle: Michaelis--Menten model with 40 observations.
Right: biofilm model with 100 observations.
}
\label{fig:dmd_sparse_all}
\end{figure}

The DMD-based augmentation strategy will be discussed here, as it was applied to the 2D linear, Michaelis--Menten, and biofilm models.~In each case, the original high-resolution trajectory was downsampled to emulate sparse experimental observations before reconstruction.~Specifically, from 500 observations, the 2D linear trajectory was reduced to 30 observations, the Michaelis--Menten trajectory to 40 observations, and the biofilm trajectory to 100 observations; these are the minimum numbers of observations needed for an accurate model to be generated.~The corresponding sparse observations are shown in Fig.~\ref{fig:dmd_sparse_all}.

The CSP diagnostics obtained from the reconstructed datasets are reported in Appendices~\ref{app:csp_2dsample}, \ref{app:csp_MMsample}, and \ref{app:csp_biofilm_sample}.~It is shown that, despite the severe reduction in temporal resolution, the reconstructed diagnostics preserve the dominant timescale structure in all cases, enabling reliable identification of the fast and slow dynamics.~The learned SIMs were compared with the corresponding analytical SIMs.~The resulting prediction errors are summarized in Table~\ref{tab:error_results_sparse_combined}.~The small errors indicate that the DMD-augmented trajectories retain sufficient dynamical information for accurate manifold identification, even when the original observations are sparse.
\begin{table}[h!]
\centering
\caption{SIM prediction errors obtained from DMD-reconstructed sparse data.}
\label{tab:error_results_sparse_combined}
\begin{tabular}{l c c c}
\toprule
\textbf{Model} & \textbf{Fast variable} & \textbf{MSE} & \textbf{RMSE} \\
\midrule
2D linear & $y$ & $3.34 \times 10^{-9}$ & $5.85 \times 10^{-5}$ \\
\hline
Michaelis--Menten, sQSSA & $c$ & $5.17 \times 10^{-5}$ & $7.19 \times 10^{-3}$ \\
Michaelis--Menten, rQSSA & $s$ & $3.43 \times 10^{-12}$ & $9.78 \times 10^{-6}$ \\
\hline
Biofilm & $G_i$ & $4.02 \times 10^{-6}$ & $2.00 \times 10^{-3}$ \\
\bottomrule
\end{tabular}
\end{table}

Table~\ref{tab:comparison_sparse_combined} compares the analytically derived reduced models with the data-driven identified models.~The identified coefficients of the models remain close to the analytical values across all systems, demonstrating that the proposed framework retains high accuracy even under sparse sampling conditions.
\begin{table}[h!]
\centering
\caption{Comparison between analytical and identified models obtained from DMD-reconstructed sparse data.}
\label{tab:comparison_sparse_combined}
\begin{tabular}{l c c}
\toprule
\textbf{Model} & \textbf{Analytical reduced model} & \textbf{Identified model} \\
\midrule
2D linear 
& $\dot{x}=-1.002x$ 
& $\dot{x}=-0.994x$ \\
\midrule
Michaelis--Menten, sQSSA 
& $\dot{s}=-0.100c$ 
& $\dot{s}=-0.098c$ \\
Michaelis--Menten, rQSSA 
& $\dot{c}=-0.100c$ 
& $\dot{c}=-0.099c$ \\
\midrule
Biofilm 
& $\dot{G}_p=1.426G_p-5.300G_pA$ 
& $\dot{G}_p=1.439G_p+\textcolor{red}{0.075G_p^2+0.065A^2}-5.554G_pA$ \\
& $\dot{A}=2.000G_p-4.000A$ 
& $\dot{A}=1.938G_p-3.895A$ \\
\bottomrule
\end{tabular}
\end{table}
Overall, it is shown that the combination of downsampling, Hankel preprocessing, and DMD-based reconstruction provides a practical mechanism for transforming sparse observations into dense trajectories suitable for data-driven multiscale model reduction.
\subsection{Preprocessing of sparse and noisy data}
\label{sec:noisydata}

In addition to sparse sampling, experimental measurements are often contaminated by noise, which can distort vector-field reconstruction, Jacobian estimation, and CSP-based diagnostics. To assess robustness under such conditions, the original trajectories are first downsampled and then contaminated with additive Gaussian noise scaled to the standard deviation of each state variable.~The noisy observations are subsequently smoothed using a third-order Savitzky–Golay filter and interpolated onto a finer temporal grid using cubic interpolation, yielding dense and smooth trajectories suitable for the proposed framework. The resulting preprocessing pipeline consists of four steps: downsampling, noise addition, Savitzky–Golay smoothing, and interpolation.

\subsubsection{Results}
\begin{figure}[h!]
\centering
\includegraphics[width=\linewidth]{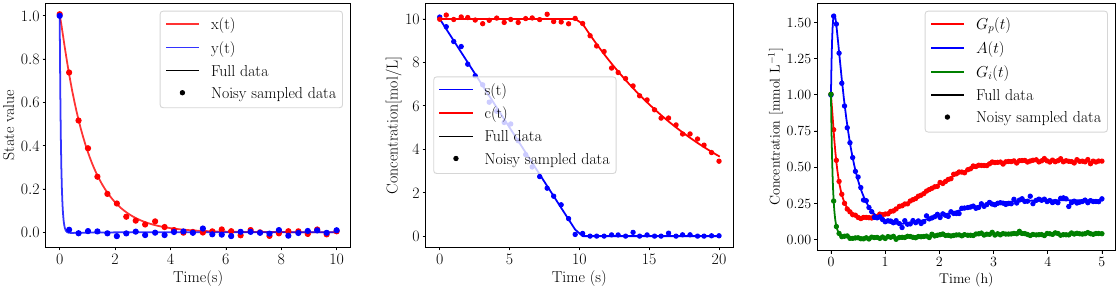}
\caption{
Full trajectories together with sparse and noisy observations used for preprocessing and interpolation.
Left: 2D linear system with 30 noisy observations.
Middle: Michaelis--Menten model with 40 noisy observations.
Right: biofilm model with 100 noisy observations.
}
\label{fig:noisy_all}
\end{figure}

The sparse and noisy preprocessing pipeline was applied to the 2D linear, Michaelis--Menten, and biofilm models. The original trajectories were downsampled to 30, 40, and 100 observations, respectively, and Gaussian noise corresponding to \(5\%\) of the standard deviation of the sampled data was added to each dataset. The resulting sparse and noisy observations are shown in Fig.~\ref{fig:noisy_all}.
The noisy observations were smoothed and interpolated before being passed to the proposed framework. The corresponding CSP diagnostics are reported in Appendices~\ref{app:csp_2dnoise}, \ref{app:csp_MMnoise}, and \ref{app:csp_biofilm_noise}. It is shown that despite the combined effects of sparse sampling and measurement noise, the dominant timescale structure was preserved in all cases, allowing the fast and slow dynamics to be identified reliably.
The learned SIMs were then compared with the corresponding analytical manifolds. The resulting prediction errors are summarized in Table~\ref{tab:error_results_noise_combined}. As expected, the errors increase relative to the noiseless sparse-data case, but the learned manifolds remain sufficiently accurate to support reduced-order model identification.
\begin{table}[h!]
\centering
\caption{SIM prediction errors obtained from sparse and noisy observations.}
\label{tab:error_results_noise_combined}
\begin{tabular}{l c c c}
\toprule
\textbf{Model} & \textbf{Fast variable} & \textbf{MSE} & \textbf{RMSE} \\
\midrule
2D linear & $y$ & $9.92 \times 10^{-5}$ & $9.95 \times 10^{-3}$ \\
Michaelis--Menten, sQSSA & $c$ & $2.39 \times 10^{-3}$ & $4.89 \times 10^{-2}$ \\
Michaelis--Menten, rQSSA & $s$ & $1.22 \times 10^{-3}$ & $3.49 \times 10^{-2}$ \\
Biofilm & $G_i$ & $1.61 \times 10^{-5}$ & $4.02 \times 10^{-3}$ \\
\bottomrule
\end{tabular}
\end{table}

Table~\ref{tab:comparison_noise_combined} compares the data-driven identified models with the analytical reduced models.~These identified models remain close to the analytical forms, although additional nonlinear terms appear in the biofilm case due to the combined effects of sparse sampling, noise, and regression over a reconstructed trajectory.~It is thus demonstrated that the proposed pre-processing and reduction framework can recover the dominant SIM structure and reduced dynamics even from sparse and noisy observations, although noise can introduce additional terms in the identified reduced equations.

\begin{table}[h!]
\centering
\caption{Comparison between analytical and identified reduced models obtained from sparse and noisy observations.}
\label{tab:comparison_noise_combined}
\begin{tabular}{l c c}
\toprule
\textbf{Model} & \textbf{Analytical reduced model} & \textbf{Identified reduced model} \\
\midrule
2D linear 
& $\dot{x}=-1.0025x$ 
& $\dot{x}=-1.135x$ \\
\midrule
Michaelis--Menten, sQSSA 
& $\dot{s}=-0.1c$ 
& $\dot{s}=-0.103c$ \\
Michaelis--Menten, rQSSA 
& $\dot{c}=-0.1c$ 
& $\dot{c}=-0.99c$ \\
\midrule
Biofilm 
& $\dot{G}_p=1.426G_p-5.3G_pA$ 
& $\dot{G}_p=1.492G_p\textcolor{red}{-0.531G_p^2}-4.523G_pA$ \\
& $\dot{A}=2G_p-4A$ 
& $\dot{A}=2.004G_p-3.979A\textcolor{red}{-0.515A^2}$ \\
\bottomrule
\end{tabular}
\end{table}



\section{Discussion}

We presented a data-driven framework in which regression on near-manifold data recovers
the slow invariant manifold (SIM) of a multiscale system, and the reduced dynamics follow by
restricting the reconstructed vector field to it. The other steps serve this one: the
surrogate vector field and its Jacobian feed CSP diagnostics, which identify the fast
variables and the time windows in which the fast modes are exhausted, and thereby select the
data the regression is fitted to. Unlike our earlier approach, which fitted separate reduced
models in each dynamical regime~\cite{muhammed2025data}, the manifold map is learned explicitly here. The output is an algebraic relation between fast
and slow variables that can be compared directly with QSSA or partial-equilibrium
approximations.

Why regression should recover the SIM has two parts. In the linear stochastic setting, the
regression slope, the Mori projection and the CSP slow eigendirection agree at leading order
in the timescale ratio $\varepsilon$, and the memory kernel reduces to its Markovian limit.
For nonlinear systems the corresponding object is Zwanzig's conditional expectation
$\mathbb{E}[y_f \mid y_s]$. In the deterministic setting of our benchmarks the argument is
simpler: once the fast modes are exhausted, trajectories lie on the SIM up to exponentially
small terms, so regression over a sufficiently rich hypothesis space is not limited to leading
order. 

The benchmarks test different aspects of this claim. The 2D linear model checks consistency
with the slow eigenspace. In the Michaelis--Menten model, the pipeline detected the switch from
the sQSSA to the rQSSA regime near $t \approx 10$\,s and the accompanying exchange of fast
variable. The biofilm model shows that the learned coefficient $\alpha = 0.06354$ matches the
higher-order SIM ($0.06350$) rather than the leading-order QSSA value $k_4/k_5 = 0.0667$. In the TMDD model, it recovered rational SIMs in two successive regimes, with two and
three fast modes, respectively. All four models are elementary mass-action networks and
therefore quadratic in the full state. It is worth noting that typically the non-polynomial expressions in the reduced models
arise from the reduction itself, which is why a polynomial library suffices for the surrogate.

Under sparse and noisy sampling the results were mixed. For the 2D, Michaelis--Menten and
biofilm models, the ordering of timescales and the CSP pointers were preserved, so fast
variables and exhaustion windows were still identified. The timescale values themselves are
bounded below by roughly the sampling interval, because data lying on the SIM carry
information about the fast rates only through $\mathcal{O}(\varepsilon)$ deviations from it.
Noise introduced spurious terms in the identified reduced equations; these originate in the
surrogate reconstruction rather than in the SIM regression, and partly compensate each other
(for the biofilm model the fixed point shifts by 1.4\%). For the TMDD model, the diagnostics
became irregular and the fast--slow split could not be identified reliably. We attribute this to the larger number of fast modes that must be resolved from an imperfect
surrogate, but a single case does not establish how data requirements scale with the
co-dimension of the fast subspace.

The framework has four main limitations. First, it assumes the full state is observed; with
unmeasured species, such as enzyme complexes, the observed dynamics are already reduced and in
general non-Markovian. Second, the accuracy of the leading-order reduction is set by the
spectral gap, with errors of order $\varepsilon = \tau_M/\tau_{M+1}$, and near-manifold data
must be available. Third, measurement noise in the slow variables biases least-squares
estimates of the manifold map, an errors-in-variables effect distinct from the intrinsic noise
treated in the stochastic analysis. Fourth, errors in the surrogate vector field propagate to
every later step. The Krylov--Arnoldi estimator needs only vector-field evaluations and scales
to large systems, but this advantage was not exercised here, since $r = N$ in all benchmarks.
The method has also not yet been applied to experimental data.

Timescale separation also underlies normative theories of how cells decide. In the
least microenvironmental uncertainty principle (LEUP) that we proposed earlier, a cell's
internal state is described by a stationary distribution conditioned on its
microenvironment, chosen to minimize the joint entropy of internal and environmental
states~\cite{hatzikirou2018leup,barua2020leup,barua2020njp}. Such a conditional
description implicitly assumes that internal states equilibrate faster than the
microenvironment changes, so the decision rule is itself a slow-manifold object: its
mean, $\mathbb{E}[y_f \mid y_s]$, is the quantity that regression on near-manifold data
estimates.  Identifying the SIM from data therefore offers a way to test this assumption and to
estimate decision maps directly from measurements. A second biomedical application is the prediction
from sparse clinical data. In the Bayesian combination of mechanistic modeling and
machine learning (BaM$^3$), the mechanistic model supplies the prior and a data-driven
regression on patient data supplies the likelihood, an approach that improved
personalized tumor growth predictions~\cite{mascheroni2021bam3}. It requires
mechanistic models simple enough to calibrate from few measurements, and the explicit
reduced models identified here are natural candidates.

.

Overall, the proposed framework provides a unified and interpretable approach for data-driven model reduction in multiscale dynamical systems. 
By unifying matrix-free local linearization, CSP-inspired multiscale diagnostics, and regression-based manifold learning, it bridges the gap between classical model reduction techniques and modern data-driven methods.
The results demonstrate that it is possible to recover both the geometric structure and reduced dynamics of complex systems directly from data, without requiring explicit knowledge of the governing equations. More broadly, the proposed framework illustrates how machine learning and mechanistic dynamical systems theory can be unified to produce physically interpretable reduced-order models, providing a foundation for future data-driven multiscale analysis of complex biological systems.

\section*{Acknowledgments}
HH would like to acknowledge the support of the RIG-2025-001. HH and DM would like to thank the UAE-NIH Collaborative Research grant AJF-NIH-25-KU.

\section*{Supporting information}

\textbf{S1 Appendix. Analytic proof of principle of the proposed framework for the 2D stochastic linear system.}\\[5pt]

\noindent\textbf{S2 Appendix. Numerical settings for the benchmark systems.}\\[5pt]

\noindent\textbf{S3 Appendix. CSP diagnostics.}

\section*{Author contributions}
\noindent\textbf{Conceptualization:} Haralampos Hatzikirou \\
\noindent\textbf{Data curation:} Dimitris M.~Manias, Ismaila Muhammed \\
\noindent\textbf{Formal analysis: } Ismaila Muhammed, Dimitris M.~Manias,  Haralampos Hatzikirou, Dimitris A. Goussis\\
\noindent\textbf{Funding acquisition:} Haralampos Hatzikirou \\
\noindent\textbf{Project administration:} Haralampos Hatzikirou \\
\noindent\textbf{Investigation:} Ismaila Muhammed, Dimitris M.~Manias,  Haralampos Hatzikirou, Dimitris A. Goussis \\ 
\noindent\textbf{Methodology:} Haralampos Hatzikirou, Dimitris A. Goussis \\
\noindent\textbf{Software:} Dimitris M.~Manias, Ismaila Muhammed \\
\noindent\textbf{Supervision:} Haralampos Hatzikirou\\
\noindent\textbf{Validation:} Dimitris A. Goussis\\
\noindent\textbf{Visualization:}  Dimitris M.~Manias, Ismaila Muhammed \\
\noindent\textbf{Writing - original draft:}  Ismaila Muhammed, Dimitris M.~Manias\\
\noindent\textbf{Writing - review \& editing: } Haralampos Hatzikirou, Dimitris A. Goussis 

\bibliography{ref}
\label{LastPageMain}

\clearpage

\setcounter{page}{1}
\renewcommand{\thepage}{\arabic{page}}

\setcounter{section}{1}
\setcounter{equation}{0}
\renewcommand{\theequation}{\arabic{section}.\arabic{equation}}

\fancyfoot[r]{\thepage/\pageref{LastPageSupp1}}

\section*{S1 Appendix. Analytic proof of principle of the proposed framework for the 2D stochastic linear system}
\label{app:covariance_analysis}

\subsection{Covariance Dynamics} 
\label{app:covariance_dynamics}
The covariance matrix
\[
\mathbf{\Sigma}(t)
=
\mathbb{E}
\left[
(\mathbf{z}_t-\boldsymbol{\mu}_t)
(\mathbf{z}_t-\boldsymbol{\mu}_t)^\top
\right]
\]
satisfies the continuous-time Lyapunov equation
\[
\dot{\mathbf{\Sigma}}
=
\mathbf{J}\mathbf{\Sigma}
+
\mathbf{\Sigma}\mathbf{J}^\top
+
\mathbf{I}.
\]

Expanding componentwise yields
\begin{align}
\dot{\Sigma}_{xx}
&=
-2a\Sigma_{xx}
+
2b\Sigma_{xy}
+
1,
\label{e:Sxx}\\
\dot{\Sigma}_{xy}
&=
c\Sigma_{xx}
-
(a+d)\Sigma_{xy}
+
b\Sigma_{yy},
\label{e:Sxy}\\
\dot{\Sigma}_{yy}
&=
2c\Sigma_{xy}
-
2d\Sigma_{yy}
+
1.
\label{e:Syy}\end{align}

At steady state,
\[
\dot{\mathbf{\Sigma}}
=
0.
\]

Setting $\dot{\mathbf{\Sigma}} = \mathbf{0}$, we extract leading-order asymptotic approximations by applying the dominant balance assumption ($d \gg a, b, c$). 

From Eq.~\eqref{e:Sxx}, assuming $\Sigma_{xy}$ is statistically small, $-2a \Sigma_{xx} + 1 \approx 0 \implies \Sigma_{xx} \approx \frac{1}{2a}$. 
From Eq.~\eqref{e:Syy}, assuming $\Sigma_{xy}$ is small relative to $d$, $-2d \Sigma_{yy} + 1 \approx 0 \implies \Sigma_{yy} \approx \frac{1}{2d}$. 
Substituting these approximations into Eq.~\eqref{e:Sxy} yields:
\begin{equation}
c \left(\frac{1}{2a}\right) - d \Sigma_{xy} + b \left(\frac{1}{2d}\right) \approx 0.
\end{equation}
Since $d \gg a$, the term $b/(2d)$ is vanishingly small compared to $c/(2a)$. The dominant balance $c \Sigma_{xx} \approx d \Sigma_{xy}$ yields $\Sigma_{xy} \approx \frac{c}{2ad}$.
Thus, the stationary covariance cloud is explicitly approximated by:
\begin{equation}
\Sigma_{xx} \approx \frac{1}{2a}, \quad \Sigma_{yy} \approx \frac{1}{2d}, \quad \text{and} \quad \Sigma_{xy} \approx \frac{c}{2ad}.
\end{equation}

\subsection{The Jacobian inverse problem}
\label{app:inverse_problem}

A standard inverse problem attempts to recover $\mathbf{J}$ from the stationary covariance $\mathbf{\Sigma}$ using the algebraic Lyapunov equation $\mathbf{J \Sigma} + \mathbf{\Sigma J^\top} + \mathbf{I} = \mathbf{0}$. Generally, this is underdetermined (providing 3 equations for 4 unknowns). However, the extreme time-scale separation naturally decouples the system, allowing us to uniquely recover the dynamically relevant features of the Jacobian.

Given a measured stationary covariance matrix, we rearrange the steady-state Lyapunov equations to isolate the parameters:
\begin{itemize}
    \item Fast rate ($d$): $d = \frac{1 + 2c \Sigma_{xy}}{2\Sigma_{yy}}\approx \frac{1}{2\Sigma_{yy}}$, since $\Sigma_{xy} \sim \mathcal{O}(1/d)$.
    \item Slow rate ($a$): $a = \frac{1 + 2b \Sigma_{xy}}{2\Sigma_{xx}}\approx \frac{1}{2\Sigma_{xx}}$, since  $2b \Sigma_{xy} \approx 0$.
    \item Forward coupling ($c$): From $c \Sigma_{xx} - (a+d)\Sigma_{xy} + b \Sigma_{yy} = 0$, we established that for $d \gg a$, we read $c \approx \frac{\Sigma_{xy}}{2\Sigma_{xx}\Sigma_{yy}}$.
\end{itemize}
The back-coupling parameter $b$ remains statistically buried in the $\mathcal{O}(1/d)$ noise floor and cannot be isolated without time-lagged data. However, $b$ is geometrically invariant to the leading-order shape of the manifold. The strong timescale separation effectively breaks the symmetry of the Lyapunov equation, allowing exact recovery of the intrinsic fast scale, slow scale, and manifold projection directly from the stationary distribution.

\subsection{Identification of the fast variable} 
\label{app:fast_variable}

We derive the right eigenvectors $\mathbf{a}_n$ by solving $\mathbf{J} \mathbf{a}_n = \lambda_n \mathbf{a}_n$. 

For the eigenvalue \(\lambda_1\) associated with the fast timescale, the corresponding right eigenvector is obtained from
\[
\begin{bmatrix}
-a-\lambda_1 & b\\
c & -d-\lambda_1
\end{bmatrix}
\begin{bmatrix}
a_1^1\\
a_1^2
\end{bmatrix}
=
\begin{bmatrix}
0\\
0
\end{bmatrix}.
\]
This implies
\begin{equation}
(-a-\lambda_1)a_1^1 + b\,a_1^2 = 0
\label{eq:first-right-eigvec}
\end{equation}
and
\begin{equation}
c\,a_1^1 + (-d-\lambda_1)a_1^2 = 0.
\label{eq:second-right-eigvec}
\end{equation}

From Eq.~\eqref{eq:first-right-eigvec}, we obtain
\begin{equation}
a_1^2=\frac{a+\lambda_1}{b}\,a_1^1.
\end{equation}
Using the leading-order approximation for \(\lambda_1\), this becomes
\begin{equation}
a_1^2=\frac{-(d-a)}{b}\,a_1^1.
\end{equation}
Since eigenvectors are defined up to a nonzero scaling factor, we choose \(a_1^1=b\), which gives
\[
a_1^2=-(d-a)\approx -d.
\]
Hence, the right eigenvector corresponding to the fast mode is
\[
\mathbf{a}_1=
\begin{bmatrix}
b\\
-d
\end{bmatrix}.
\]

The left eigenvector is obtained from
\[
\begin{bmatrix}
-a-\lambda_1 & c\\
b & -d-\lambda_1
\end{bmatrix}
\begin{bmatrix}
b_1^1\\
b_2^1
\end{bmatrix}
=
\begin{bmatrix}
0\\
0
\end{bmatrix},
\]
which is equivalent to solving
\[
\mathbf{J}^T(\mathbf{b}^{1})^T=\lambda_1(\mathbf{b}^{1})^T.
\]
This implies
\begin{equation}
(-a-\lambda_1)b_1^1 + c\,b_2^1 = 0
\label{eq:first-left-eigvec}
\end{equation}
and
\begin{equation}
b\,b_1^1 + (-d-\lambda_1)b_2^1 = 0.
\label{eq:second-left-eigvec}
\end{equation}

From Eq.~\eqref{eq:first-left-eigvec}, we obtain
\begin{equation}
b_2^1=\frac{a+\lambda_1}{c}\,b_1^1.
\end{equation}
Using the leading-order approximation for \(\lambda_1\), this becomes
\begin{equation}
b_2^1=\frac{-(d-a)}{c}\,b_1^1.
\end{equation}
Choosing \(b_1^1=c\), we get
\[
b_2^1=-(d-a)\approx -d.
\]
Hence, the (unnormalized) left eigenvector corresponding to the fast mode points along
\[
\mathbf{b}^{1}\propto
\begin{bmatrix}
c & -d
\end{bmatrix}.
\]

To satisfy the biorthogonality condition $\mathbf{b}^{1}\mathbf{a}_{1} = 1$, we must jointly normalize the eigenvectors. Given our choice of right eigenvector $\mathbf{a}_1 \approx \begin{bmatrix} b & -d \end{bmatrix}^\top$, the appropriately normalized left eigenvector is:
\begin{equation}
\mathbf{b}^{1} \approx \frac{1}{d^2}\begin{bmatrix}c & -d\end{bmatrix}.
\end{equation}
This ensures $\mathbf{b}^1\mathbf{a}_1 = \frac{1}{d^2}(cb + d^2) = \frac{cb}{d^2} + 1 \approx 1$ asymptotically.

The CSP pointer for the fast mode identifies which physical variable drives the fast dynamics. Given that 
\[
\mathbf{a}_{1}\mathbf{b}^{1}=
\begin{bmatrix}
b\\
-d
\end{bmatrix}
\frac{1}{d^2}
\begin{bmatrix}
c & -d
\end{bmatrix}
=
\begin{bmatrix}
bc/d^2 & -b/d\\
-c/d & 1
\end{bmatrix},
\]
the pointer associated with the fast mode is calculated as:
\begin{equation}
Po^1=\operatorname{diag}(\mathbf{a}_{1}\mathbf{b}^{1})
=
\begin{bmatrix}
bc/d^2\\
1
\end{bmatrix}
\approx
\begin{bmatrix}
0\\
1
\end{bmatrix},
\end{equation}

which indicates that the second variable is the fast variable, since the pointer associated with the fast mode is fully aligned with it.

\subsection{The Mori--Zwanzig projection of the 2D linear system}
\label{app:MZ}

This appendix explains, from first principles, why the regression-based manifold and the reduced slow dynamics identified by the framework are not coincidences but exact consequences of the classical \emph{Mori--Zwanzig} (MZ) projection \cite{mori1965transport,zwanzig1961memory,zwanzig1973nonlinear,chorin2000optimal}. No prior familiarity with MZ theory is assumed: we build the necessary vocabulary on the simplest possible example---the very 2D linear system of the main text---and then show that three objects that look different (a CSP eigenvector, a least-squares regression line, and a statistical projection) are the same.

Throughout we use the system in component form, with unit-intensity white noise $\eta=(\eta_x,\eta_y)$:
\begin{equation}
\dot{x}=-a\,x+b\,y+\eta_x,
\qquad
\dot{y}=c\,x-d\,y+\eta_y,
\qquad a,b,c,d>0,\quad d\gg a,b,c.
\label{eq:MZ_system}
\end{equation}

\subsubsection{The problem, in plain terms}

We have two variables. One of them, $y$, is \emph{fast}: it has the large damping rate $d$, so it forgets its own past almost instantly and is always being pulled toward whatever value the slow variable $x$ currently dictates. The other, $x$, is \emph{slow}. We would like to \emph{eliminate} $y$ and keep a single equation for $x$ alone.

The naive way is to set $\dot{y}=0$ (``$y$ has reached steady state''), solve $cx-dy=0$ for $y=\tfrac{c}{d}x$, and substitute. This is exactly the quasi-steady-state (QSS) recipe and it gives the correct \emph{deterministic} answer, but it hides two things. First, $y$ does not respond to $x$ \emph{instantaneously}; it responds with a short delay of order $1/d$, so the true equation for $x$ depends not only on $x(t)$ but on a little of its recent \emph{past} (a \emph{memory} term). Second, $y$ is itself kicked by the noise $\eta_y$, and that randomness leaks into the $x$-equation as an extra \emph{fluctuating force} that $x$ cannot predict. The MZ formalism is the bookkeeping that makes these two hidden pieces explicit and exact; when $d$ is large (our regime), both shrink in a controlled way and we recover the clean reduced model, now \emph{knowing} what was neglected and how small it is.

\subsubsection{Observables as vectors; covariance as an inner product}

To speak of ``the best estimate of $y$ given $x$'' we need a notion of \emph{length} and \emph{angle} between random variables, and the statistics of the system supply one for free. For two scalar observables $A$ and $B$ (here, linear combinations of $x$ and $y$), define their inner product as the stationary correlation
\begin{equation}
\langle A,\,B\rangle:=\mathbb{E}_{\mathrm{ss}}[A\,B],
\label{eq:MZ_inner}
\end{equation}
where $\mathbb{E}_{\mathrm{ss}}$ is the average over the long-time (stationary) distribution. Because the means decay to zero (Sec.~\ref{app:covariance_dynamics}, S1), this is just the covariance; in particular $\langle x,x\rangle=\Sigma_{xx}$, $\langle y,y\rangle=\Sigma_{yy}$, and $\langle x,y\rangle=\Sigma_{xy}$. This $\langle\cdot,\cdot\rangle$ behaves like an ordinary dot product: it is symmetric, linear, and $\langle A,A\rangle=\mathbb{E}[A^2]\ge 0$ defines a length $\|A\|=\sqrt{\langle A,A\rangle}$. Two observables are \emph{orthogonal} when $\langle A,B\rangle=0$, i.e. when they are \emph{uncorrelated}. In short: random variables become geometric vectors once ``how aligned are $A$ and $B$'' means ``how correlated are they''; the covariance matrix is the metric.

\subsubsection{Projection = best linear regression = conditional mean}

The \emph{Mori projection} of any observable onto the slow variable $x$ is the orthogonal projection in the inner product \eqref{eq:MZ_inner}:
\begin{equation}
\mathcal{P}A:=\frac{\langle A,\,x\rangle}{\langle x,\,x\rangle}\,x,
\qquad
\mathcal{Q}:=\mathcal{I}-\mathcal{P}.
\label{eq:MZ_proj}
\end{equation}
$\mathcal{P}A$ is the multiple of $x$ closest to $A$ in the norm $\|\cdot\|$. The leftover $\mathcal{Q}A=A-\mathcal{P}A$ is by construction \emph{uncorrelated with $x$}, since $\langle A-\tfrac{\langle A,x\rangle}{\langle x,x\rangle}x,\,x\rangle=\langle A,x\rangle-\langle A,x\rangle=0$. Three facts make $\mathcal{P}$ exactly the operation a regression performs. It is least squares: $\mathcal{P}A$ minimizes $\mathbb{E}[(A-\beta x)^2]$ over slopes $\beta$, with minimizer $\beta=\langle A,x\rangle/\langle x,x\rangle=\Sigma_{Ax}/\Sigma_{xx}$, the textbook regression coefficient. It is the conditional mean: for jointly Gaussian zero-mean variables $\mathbb{E}[A\mid x]=\tfrac{\Sigma_{Ax}}{\Sigma_{xx}}x$, identical to \eqref{eq:MZ_proj}, and the stationary state of the linear system \eqref{eq:MZ_system} is jointly Gaussian. And it is idempotent, $\mathcal{P}^2=\mathcal{P}$, as any projection must be. In words: $\mathcal{P}y$ is ``the part of the fast variable that the slow variable can explain'' and $\mathcal{Q}y$ is ``the part it cannot''; fitting a line $y\approx\beta x$ to data \emph{is} computing $\mathcal{P}y$.

\subsubsection{Exact elimination of the fast variable}

We eliminate $y$ from \eqref{eq:MZ_system} \emph{without approximation}. The $y$-equation is a linear first-order ODE, $\dot{y}+d\,y=c\,x+\eta_y$. Multiplying by the integrating factor $e^{dt}$ makes the left side a total derivative, $\tfrac{d}{dt}(e^{dt}y)=e^{dt}(cx+\eta_y)$; integrating from $0$ to $t$ and dividing by $e^{dt}$ gives the exact solution for the fast variable as a functional of the slow history,
\begin{equation}
y(t)=\underbrace{e^{-dt}y(0)}_{\text{memory of initial }y}
+\int_0^t e^{-d(t-s)}\big(c\,x(s)+\eta_y(s)\big)\,ds.
\label{eq:MZ_ysol}
\end{equation}
Equation \eqref{eq:MZ_ysol} says $y(t)$ is a \emph{weighted average of the past}, the weight $e^{-d(t-s)}$ emphasizing the immediate past (within a window of width $1/d$) and exponentially forgetting older times; this window width $1/d$ is the fast time scale $\tau_1$. Substituting \eqref{eq:MZ_ysol} into the $x$-equation $\dot{x}=-ax+by+\eta_x$ yields the closed, exact \emph{generalized Langevin equation} (GLE)
\begin{equation}
\dot{x}(t)=-a\,x(t)+bc\!\int_0^t e^{-d(t-s)}x(s)\,ds+F(t),
\label{eq:MZ_GLE}
\end{equation}
with everything independent of the slow history collected into
\begin{equation}
F(t)=\underbrace{b\,e^{-dt}y(0)}_{\text{initial-layer transient}}
+\underbrace{b\!\int_0^t e^{-d(t-s)}\eta_y(s)\,ds+\eta_x(t)}_{\text{filtered fast noise + direct noise}}.
\label{eq:MZ_force}
\end{equation}
No information has been lost; the dynamics has only been reorganized.


\subsubsection{Anatomy of the generalized Langevin equation}

The GLE~(1.20) decomposes the exact evolution of $x(t)$ into three distinct physical contributions. The bare streaming term $-a\,x(t)$ is the instantaneous (Markovian) drift---how $x$ would evolve in the absence of coupling to $y$. The memory term $\int_0^t K(t-s)x(s)\,ds$, with exponential memory kernel
\begin{equation}
K(t) = bc\,e^{-dt},
\label{eq:S1_kernel}
\end{equation}
captures the delayed feedback that the slow variable exerts on itself through the fast variable: $x$ drives $y$ at rate $c$, $y$ feeds back into $x$ at rate $b$, and the round-trip relaxation occurs over the fast time scale $\tau_1 = 1/d$. 

Finally, the forcing term $F(t)$ in~(1.21) collects the contributions not generated by the history of $x$: the initial-layer transient $b e^{-dt} y(0)$, the exponentially filtered fast noise $b \int_0^t e^{-d(t-s)}\eta_y(s)\,ds$, and the direct slow white noise $\eta_x(t)$. Note that because $F(t)$ is obtained here by direct variation of parameters rather than a formal Mori orthogonal dynamics propagator, $F(0) = b y(0) + \eta_x(0)$ is not strictly orthogonal to $x(0)$ in the stationary inner product; specifically, decomposing $y(0) = \mathcal{P}y(0) + \mathcal{Q}y(0)$ yields a correlated transient piece $b e^{-dt}\mathcal{P}y(0) = (bc/d)e^{-dt}x(0)$ alongside the genuinely orthogonal initial fluctuation $b e^{-dt}\mathcal{Q}y(0)$. However, this correlated transient decays exponentially at the fast rate $d$ and vanishes identically if the trajectory is initialized on the slow manifold ($y(0) = (c/d)x(0)$ absorbed into the lower limit of the memory integral). Beyond the initial fast boundary layer ($t \gg 1/d$), $F(t)$ reduces purely to the unpredictable stochastic forcing.

\subsubsection{Why memory and noise collapse when the fast rate is large}

Here the stiffness $d\gg a,b,c$ earns its keep. In the memory integral the kernel only weights a window of width $1/d$, over which a slow variable barely changes, so we may evaluate $x(s)$ at $s=t$:
\begin{equation}
\int_0^t e^{-d(t-s)}x(s)\,ds\approx x(t)\int_0^t e^{-d(t-s)}\,ds
=x(t)\,\frac{1-e^{-dt}}{d}\;\xrightarrow[\;t\gg 1/d\;]{}\;\frac{x(t)}{d}.
\label{eq:MZ_markov}
\end{equation}
The non-local memory term becomes a local term proportional to $x(t)$ with coefficient $bc/d$: the kernel has effectively become a delta function, $K(t)\to(\int_0^\infty K)\,\delta(t)$ with weight $\int_0^\infty K(t)\,dt=bc/d$. Simultaneously the fluctuating force \eqref{eq:MZ_force} decays at rate $d$, so after the initial fast layer it has relaxed to a small, rapidly decorrelating residual. Dropping it and using \eqref{eq:MZ_markov}, the deterministic skeleton of \eqref{eq:MZ_GLE} becomes
\begin{equation}
\dot{x}=\Big(-a+\frac{bc}{d}\Big)x=\lambda_2\,x,
\label{eq:MZ_reduced}
\end{equation}
precisely the reduced slow dynamics of the main text: the first term $-a$ is the bare slow rate and $+bc/d$ is the \emph{integrated memory}, the small systematic feedback the fast variable contributes before it is slaved. Thus ``memory collapses to Markovian'' and ``the fast CSP mode is exhausted ($f^1\approx 0$)'' are the same statement in two dialects: after a time $\sim 1/d$ the fast variable carries no independent information, being fully determined by $x$, and its only lasting effect is a tiny renormalization of the slow rate.

\subsubsection{The slow eigenvalue is the Mori frequency}

A second, purely statistical route to the same coefficient does not pass through the kernel integral. The MZ \emph{streaming frequency} is the projection of the slow velocity onto the slow variable, $\Omega=\langle\dot{x},x\rangle/\langle x,x\rangle$. Since the white noise is uncorrelated with the present state, $\langle\eta_x,x\rangle=0$, the numerator is $\langle\dot{x},x\rangle=-a\Sigma_{xx}+b\Sigma_{xy}$, so that, using $\Sigma_{xy}/\Sigma_{xx}=c/d$,
\begin{equation}
\Omega=-a+b\,\frac{\Sigma_{xy}}{\Sigma_{xx}}=-a+\frac{bc}{d}=\lambda_2
=\underbrace{(-a)}_{\text{bare streaming}}+\underbrace{\int_0^\infty K(t)\,dt}_{\text{integrated memory}}.
\label{eq:MZ_freq}
\end{equation}
The two derivations agree: the effective slow rate is the bare rate dressed by the memory, and it equals the slow eigenvalue $\lambda_2$.

\subsubsection{The slow invariant manifold \emph{is} the Mori projection}

Recall the stationary covariances (Sec.~\ref{app:covariance_dynamics} S1), $\Sigma_{xx}\approx\tfrac{1}{2a}$, $\Sigma_{yy}\approx\tfrac{1}{2d}$, $\Sigma_{xy}\approx\tfrac{c}{2ad}$. Applying the projection \eqref{eq:MZ_proj} to the fast variable $y$,
\begin{equation}
\mathcal{P}y=\frac{\langle y,x\rangle}{\langle x,x\rangle}\,x
=\frac{\Sigma_{xy}}{\Sigma_{xx}}\,x
=\frac{c/(2ad)}{1/(2a)}\,x
=\frac{c}{d}\,x
=\mathbb{E}[y\mid x].
\label{eq:MZ_Py}
\end{equation}
Now compare with what the framework's \emph{regression} computes. Fitting $y=\beta_1 x+\beta_0+\epsilon$, the intercept vanishes because the means are zero, $\beta_0=\mu_y-\beta_1\mu_x=0$, and the slope is the standard estimator
\begin{equation}
\beta_1=\frac{\mathrm{Cov}(x,y)}{\mathrm{Var}(x)}=\frac{\Sigma_{xy}}{\Sigma_{xx}}\approx\frac{c}{d}.
\label{eq:MZ_beta}
\end{equation}
Equations \eqref{eq:MZ_Py} and \eqref{eq:MZ_beta} are identical: \emph{the learned SIM $y=\tfrac{c}{d}x$ is the Mori projection $\mathcal{P}y$}. It also agrees with the geometric/CSP answer, since the exhausted-mode condition $\mathbf{b}^{1}\dot{\mathbf{z}}=\dot{y}=cx-dy\approx 0$ gives the same line and the slow right eigenvector is $\mathbf{a}_2\propto(1,\,c/d)^\top$. Three constructions, one manifold.

\subsubsection{What the regression discards lives off the manifold}

The fluctuation the regression \emph{removes} is exactly the orthogonal part $\mathcal{Q}y=y-\tfrac{c}{d}x$. It is uncorrelated with the slow variable,
\[
\langle\mathcal{Q}y,\,x\rangle=\Sigma_{xy}-\frac{\Sigma_{xy}}{\Sigma_{xx}}\Sigma_{xx}=0,
\]
so no information about the manifold is discarded---only the off-manifold scatter, whose size is the residual variance
\begin{equation}
\langle\mathcal{Q}y,\,\mathcal{Q}y\rangle=\Sigma_{yy}-\frac{\Sigma_{xy}^2}{\Sigma_{xx}}
=\frac{1}{2d}\Big(1-\frac{c^2}{ad}\Big)\approx\frac{1}{2d}.
\label{eq:MZ_residual}
\end{equation}
This is the precise content of the statement that ``the noise cloud aligns its principal axis with the slow subspace'': the regression returns the conditional mean---the \emph{minimum-variance} estimate---which by definition is blind to the perpendicular fast scatter \eqref{eq:MZ_residual}. Since that scatter is $\mathcal{O}(d^{-1})$, the manifold is recovered with an error that shrinks as the system becomes stiffer. Geometrically: the stationary data form an elongated ellipse in the $(x,y)$ plane whose long axis lies almost along $y=\tfrac{c}{d}x$ and whose short axis (thickness $\sim 1/\sqrt{2d}$) is the fast noise; least squares fits the long axis and ignores the thickness, which is why fitting a line to noisy fast--slow data automatically denoises.

\subsubsection{Shared fast decay scale across memory, noise, and variance}

Although the system is generally non-reversible ($b \neq c$, so detailed balance and the classical thermodynamic second fluctuation--dissipation theorem do not strictly apply), the memory kernel and the colored component of the fluctuating force are governed by the exact same fast relaxation rate $d$. Specifically, while the memory kernel decays as $K(t) = bc\,e^{-dt}$, the stationary autocorrelation of the filtered fast noise $F_y(t) := b\int_{-\infty}^t e^{-d(t-s)}\eta_y(s)\,ds$ evaluates to
\begin{equation}
\langle F_y(t) F_y(0) \rangle = \frac{b^2}{2d}\,e^{-d|t|}.
\label{eq:S1_force_corr}
\end{equation}
Thus, even though the prefactors $bc$ and $b^2/(2d)$ differ due to the asymmetric coupling, the colored noise decorrelates on the exact same fast time scale $\tau_1 = 1/d$ over which the memory kernel acts, because both originate from eliminating the same fast linear mode $y$. At the static level, the algebraic Lyapunov equation $\mathbf{J}\mathbf{\Sigma} + \mathbf{\Sigma} \mathbf{J}^\top + \mathbf{I} = \mathbf{0}$ encodes the corresponding stationary covariance balance that fixes the inner product and the projection $\mathcal{P}$. A single fast parameter, $d$, therefore simultaneously controls the memory decay time $\tau_{\text{mem}} = 1/d$, the colored-noise correlation time $1/d$, and the off-manifold variance $\langle \mathcal{Q}y, \mathcal{Q}y \rangle \approx 1/(2d)$.
\subsubsection{Regression slope versus SIM slope convergence} 
 
\label{SI:convergence}

Let $\delta := a/d\ll 1$, with $a, b, c$ fixed, $c \neq 0$, and $d \to \infty$. As in
Section~1.3, we chose the slow right eigenvector $\mathbf{a}_2 = (1,m)^\top$ satisfies
$\mathbf{J}\mathbf{a}_2 = \lambda_2\mathbf{a}_2$, so that the
SIM is $y = m\, x$. We show that
\begin{equation*}
  \frac{\beta_1}{m} = 1 + \delta\left(\frac{b}{c} - 2\right) + \mathcal{O}(\delta^2).
\end{equation*}
Since both slopes are of order $\delta$, each must be expanded to second order.
 
\paragraph{Regression slope.} At steady state, Eq~(1.1) divided by $\Sigma_{xx}$ gives
$1/\Sigma_{xx} = 2(a - b\beta_1)$, and Eq~(1.3) gives
$\Sigma_{yy} = (2c\,\Sigma_{xy} + 1)/(2d)$. Substituting both into Eq~(1.2), dividing by
$\Sigma_{xx}$ and multiplying by $d$ yields the exact result
\begin{equation*}
  \beta_1 = \frac{cd + ab}{d^2 + ad - bc + b^2}.
\end{equation*}
Dividing numerator and denominator by $d^2$ and using $1/d = \delta/a$,
\begin{equation}
  \beta_1 = \frac{\delta}{a}\,\frac{c + \delta b}{1 + \delta + \mathcal{O}(\delta^2)}
          = \frac{\delta}{a}\,\bigl[\,c + \delta\,(b - c) + \mathcal{O}(\delta^2)\bigr].
\end{equation}
 
\paragraph{SIM slope.} The two components of
$\mathbf{J}\mathbf{a}_2 = \lambda_2\mathbf{a}_2$ read
$-a - \lambda_2 + b\,m = 0$ and $c - (d + \lambda_2)\,m = 0$, so
$\lambda_2 = -a + b\,m$ and $m\,(d - a + b\,m) = c$. Dividing by $d$,
\begin{equation*}
  m\,\Bigl(1 - \delta + \frac{b\,m}{d}\Bigr) = \frac{\delta}{a}\,c .
\end{equation*}
Since $m= \mathcal{O}(\delta)$, the term $b\,m/d$ is $\mathcal{O}(\delta^2)$, and
\begin{equation}
  m= \frac{\delta}{a}\,\bigl[\,c\,(1 + \delta) + \mathcal{O}(\delta^2)\bigr].
\end{equation}
 
\paragraph{Ratio.} Dividing 
(1.30) 
by 
(1.31) 
cancels the common factor $\varepsilon/a$:
\begin{align*}
  \frac{\beta_1}{m}
  &= \frac{c + \delta\,(b - c) + \mathcal{O}(\delta^2)}{c\,(1 + \delta) + \mathcal{O}(\delta^2)}
   = \Bigl[1 + \delta\Bigl(\frac{b}{c} - 1\Bigr)\Bigr](1 - \delta) + \mathcal{O}(\delta^2)\\
  &= 1 + \delta\Bigl(\frac{b}{c} - 2\Bigr) + \mathcal{O}(\delta^2).
\end{align*}

\paragraph{Remarks.}
(i) The term $b/c$ comes from $b\,\Sigma_{yy}$ in Eq~(1.2), the term neglected in the
dominant balance of Eq~(1.4). It is fast noise transmitted to $x$ through the
back-coupling and has no counterpart in the SIM at this order; the term $-2$ collects the
$\mathcal{O}(\delta)$ corrections to the two slopes, $-\delta c$ in 
(1.30) 
and
$+\delta c$ in 
(1.31)
.
(ii) The $\mathcal{O}(\delta^2)$ remainder also contains terms in $b/c$, so the
expansion is informative only when $\delta\,|b/c| \ll 1$. For $c = 0$ the ratio is
undefined: $y = 0$ is invariant and $m = 0$, whereas $\beta_1 = ab/(d^2 + ad + b^2)
\neq 0$, a slope produced by the noise alone.
(iii) If the identity $\mathbf{I}$ in the Lyapunov equation is replaced by a diagonal noise
covariance $\operatorname{diag}(\sigma_x^2, \sigma_y^2)$, the same steps give
$b\,\sigma_y^2/\sigma_x^2$ in place of $b$. Regression and SIM then agree at leading
order only if $\sigma_y^2/\sigma_x^2$ stays bounded as $d \to \infty$.

\subsubsection{Dictionary and scope}

Table~\ref{tab:MZ_dictionary} makes the equivalence explicit, so a reader fluent in either language can translate.

\begin{table}[H]
\centering
\caption{Correspondence between the CSP/framework language and the Mori--Zwanzig language for the 2D linear system.}
\label{tab:MZ_dictionary}
\small
\setlength{\tabcolsep}{4pt}
\begin{tabular}{@{}p{2.4cm}p{3.2cm}p{4.0cm}p{2.5cm}@{}}
\toprule
Concept & CSP / framework & Mori--Zwanzig & This system \\
\midrule
Variable split & fast vs.\ slow variable & unresolved $\mathcal{Q}$ vs.\ resolved $\mathcal{P}$ & $y$ vs.\ $x$ \\
Fast-variable ID & CSP pointer $Po^1$ on $y$ & localization of $\mathrm{range}(\mathcal{Q})$ & $Po^1\!\approx\!(0,1)^\top$ \\
Manifold & regression $y=\beta_1 x$ & Mori projection $\mathcal{P}y=\mathbb{E}[y\mid x]$ & $y=\tfrac{c}{d}x$ \\
Off-manifold scatter & residual $\epsilon$ & $\mathcal{Q}y$ (fluctuating force) & variance $\approx\tfrac{1}{2d}$ \\
Reduced rate & slow eigenvalue $\lambda_2$ & streaming freq.\ $\Omega=-a+\int_0^\infty\!K$ & $-a+\tfrac{bc}{d}$ \\
Mode exhausted & $f^1\approx 0$ & memory kernel collapses & window $\sim 1/d$ \\
Statistics fixer & Lyapunov equation & stationary covariance balance (shared decay rate d) & Eqs.~\eqref{eq:MZ_residual},\,\eqref{eq:S1_force_corr} \\
\bottomrule
\end{tabular}
\end{table}

A few caveats delimit the clean picture. The \emph{exact} objects are the elimination \eqref{eq:MZ_ysol}, the GLE \eqref{eq:MZ_GLE}, and the kernel $K(t)=bc\,e^{-dt}$: these hold for any $a,b,c,d>0$, with no approximation. The \emph{leading-order} results are the covariance estimates, the reduced model \eqref{eq:MZ_reduced}, and the coincidence of the projection direction with the slow eigenvector; their relative error is $\mathcal{O}(bc/ad)$, vanishing as $d\to\infty$. More broadly, for \emph{linear, Gaussian} systems the Mori projection coincides with the conditional mean and is computed exactly by linear regression, which is why this benchmark is so transparent. For \emph{nonlinear} systems the optimal projection is still the conditional expectation $\mathbb{E}[y_f\mid y_s]$, but it is generally a nonlinear function and the memory kernel is no longer a simple exponential; this is exactly the setting the main framework addresses empirically, by replacing linear regression with symbolic/nonlinear regression for the manifold map $h^\star$ and by using CSP diagnostics to decide \emph{when} the memory has effectively collapsed. The linear analysis here is therefore the rigorous, fully solvable limit that explains the mechanism, while the nonlinear benchmarks of the Applications section are its practical generalization.

\label{LastPageSupp1}

\clearpage

\setcounter{page}{1}
\renewcommand{\thepage}{\arabic{page}}
\setcounter{section}{2}
\setcounter{subsection}{0}
\setcounter{equation}{0}
\renewcommand{\theequation}{\arabic{section}.\arabic{equation}}

\fancyfoot[r]{\thepage/\pageref{LastPageSupp2}}

\section*{S2 Appendix. Numerical settings for the benchmark systems}
\label{app:numerical_settings}

This appendix summarizes the numerical settings, parameter values, initial conditions, and sampling procedures used to generate the reference trajectories for all benchmark systems considered in the Applications section.

\subsection{2D linear model}

The reference trajectories for the illustrative linear system were generated from
\[
\dot x = -x + 0.5\,y,
\qquad
\dot y = -0.1\,x - 20\,y.
\]
The system was integrated over
\[
t\in[0,10],
\]
using an adaptive LSODA-type solver with relative and absolute tolerances
\[
\texttt{rtol}=10^{-9},
\qquad
\texttt{atol}=10^{-12}.
\]
The numerical solution was sampled uniformly at
\[
T=500
\]
time instances. The initial condition was chosen as
\[
\mathbf{z}(0)=(1.0,1.0)^\top.
\]


\subsection{Michaelis--Menten model}
The governing equations are
\[
\begin{aligned}
\dot{s} &= -k_{1f}(e_0-c)s + k_{1b}c,\\
\dot{c} &= k_{1f}(e_0-c)s-(k_{1b}+k_2)c.
\end{aligned}
\]
The parameter values were set to
\[
\begin{aligned}
k_{1f}&=10000,
&
k_{1b}&=0.01,
&
k_2&=0.1,
&
e_0&=10.
\end{aligned}
\]
The system was integrated over
\[
t\in[0,20],
\]
using a high-accuracy stiff solver with
\[
\texttt{rtol}=10^{-12},
\qquad
\texttt{atol}=10^{-18}.
\]
The initial conditions were
\[
s(0)=10,
\qquad
c(0)=9.999999.
\]

\subsection{Biofilm model}
\label{app:numerical_settings_Biofilm}

The biofilm model governing equations are
\[
\begin{aligned}
\dot G_p &= k_1\,G_E\,G_p - k_4\,G_p - k_2\,A\,G_p,\\
\dot A   &= -k_3\,A + k_5\,G_i,\\
\dot G_i &= k_4\,G_p - k_5\,G_i.
\end{aligned}
\]
The parameter values were set to
\[
\begin{aligned}
k_1 &= 0.3426,
&
k_2 &= 5.3,
&
k_3 &= 4.0,\\
k_4 &= 2.0,
&
k_5 &= 30.0,
&
G_E &= 10.0.
\end{aligned}
\]
The system was integrated over
\[
t\in[0,5]\;\text{hours},
\]
using an adaptive LSODA-type stiff solver with
\[
\texttt{rtol}=10^{-9},
\qquad
\texttt{atol}=10^{-12}.
\]
The numerical solution was sampled uniformly at
\[
T=1000
\]
time instances.

The initial condition was
\[
\mathbf{y}(0)=
\big(G_p(0),A(0),G_i(0)\big)^\top
=
(1.0,1.0,1.0)^\top.
\]


\subsection{Target-mediated drug disposition model}

The TMDD model governing equations are

\[
\begin{aligned}
\dot L_d &= -k_a L_d,\\
\dot L   &= k_a L_d - k_{el} L - k_{pt} L 
            - k_{on} L\,R + k_{off} RL + k_{tp} L_t,\\
\dot L_t &= k_{pt} L - k_{tp} L_t,\\
\dot R   &= -k_{on} L\,R + k_{off} RL + k_{syn} - k_{deg} R,\\
\dot{RL} &= k_{on} L\,R - k_{off} RL - k_{int} RL.
\end{aligned}
\]

The parameter values were set to

\[
\begin{aligned}
k_a &= 0.5,
&
k_{el} &= 0.05,
&
k_{pt} &= 0.1,\\
k_{on} &= 100,
&
k_{off} &= 10,
&
k_{tp} &= 0.1,\\
k_{syn} &= 2,
&
k_{deg} &= 5,
&
k_{int} &= 10.
\end{aligned}
\]

The system was integrated over
\[
t\in[0,50],
\]
using an adaptive stiff solver with strict relative and absolute tolerances.

The numerical solution was sampled uniformly at
\[
T=500
\]
time instances.

The initial conditions were chosen as
\[
\mathbf{y}(0)=
\big(L_d(0),L(0),L_t(0),R(0),RL(0)\big)^\top
=
\left(
\frac{D_1}{V},
\frac{D_2}{V},
0,
\frac{k_{syn}}{k_{deg}},
0
\right)^\top,
\]
with
\[
D_1=100,
\qquad
D_2=0,
\qquad
V=3.
\]


\label{LastPageSupp2}

\clearpage
\setcounter{page}{1}
\renewcommand{\thepage}{\arabic{page}}
\setcounter{section}{3}
\setcounter{subsection}{0}
\setcounter{subsubsection}{0}

\setcounter{equation}{0}
\renewcommand{\theequation}{\arabic{section}.\arabic{equation}}

\fancyfoot[r]{\thepage/\pageref{LastPageSupp3}}

\section*{S3 Appendix. CSP diagnostics}
\label{app:csp_diagnostics}

\subsection{2D linear model}
\label{app:csp_2d}

\subsubsection{Reference diagnostics}
\label{app:csp_2d_reference_diagnostics}
\begin{figure}[H]
\centering
\includegraphics[width=0.88\linewidth]{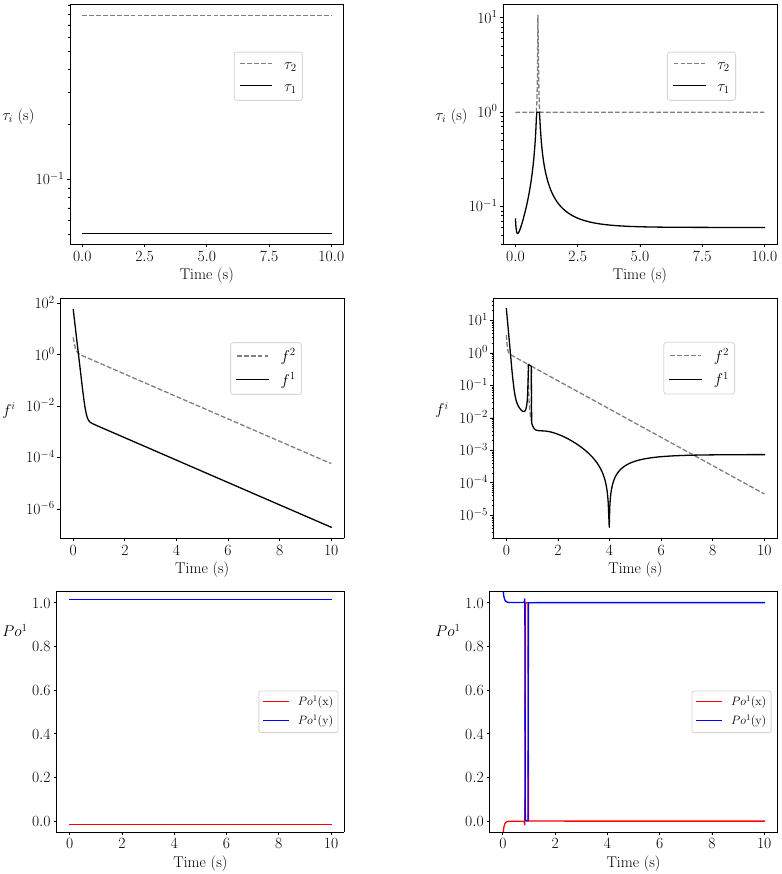}
\caption{
Comparison of the developing timescales (top), CSP amplitudes (middle), and pointers (bottom) for the 2D linear model, as obtained analytically (left) and through the proposed framework (right).~Black line refers to the fast exhausted mode.}
\label{fig:csp_analysis_2D}
\end{figure}

Figure \ref{fig:csp_analysis_2D}~(top) presents the developing time scales of the 2D linear model as calculated by both the analytical model (left) and the proposed framework (right). A time scale gap is evident in both, indicative of the multi-scale character of the model. It is shown that, in both setups, the amplitude of the fastest dissipative mode becomes negligible (middle), indicating the fast mode exhaustion and forming equilibria, based on which the reduced model can be constructed.~The CSP  pointer identifies variable $x$ to be associated directly to the fast exhaused mode; $x$ is a fast variable and a valid QSSA can be assumed for $x$.

\FloatBarrier

\subsubsection{Downsampled diagnostics}
\label{app:csp_2dsample}
The diagnostics obtained on downsampled data by the proposed framework are not affected qualitatively as compared to the ones obtained by the analytical model.~A valid QSSA can be assumed for $x$.


\begin{figure}[H]
\centering
\includegraphics[width=0.88\linewidth]{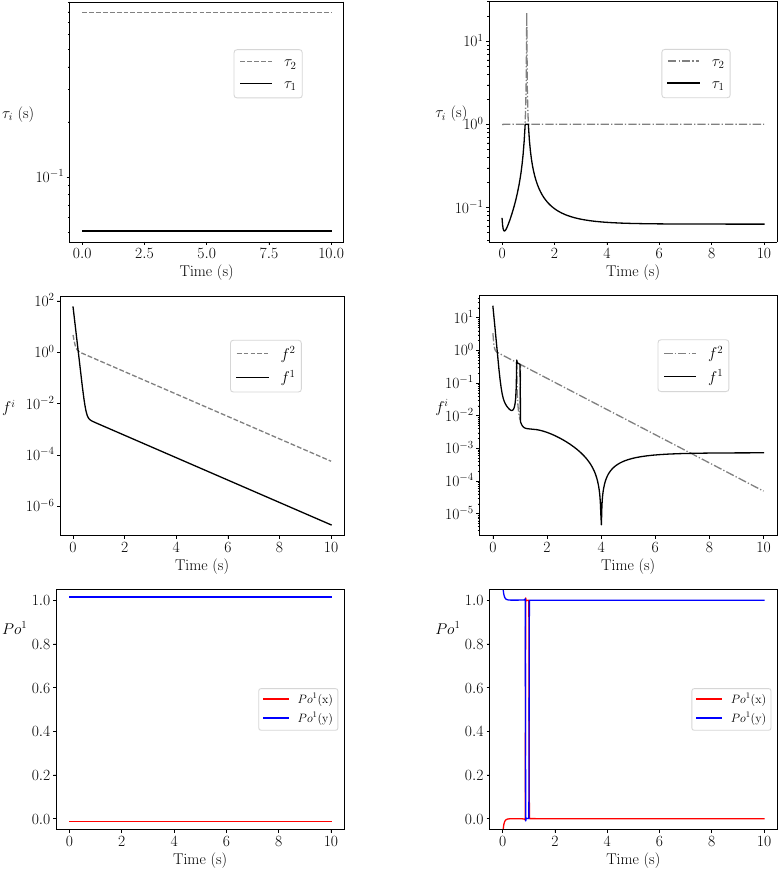}
\caption{
Comparison of the developing timescales (top), CSP amplitudes (middle), and pointers (bottom) for the 2D linear model, as obtained analytically (left) and through the proposed framework (right).~Black line refers to the fast exhausted mode.}
\label{fig:csp_analysis_2Dsample}
\end{figure}

\subsubsection{Downsampled \& noisy diagnostics}
\label{app:csp_2dnoise}
The diagnostics obtained on downsampled and noisy data by the proposed framework are affected slightly qualitatively as compared to the ones obtained by the analytical model (left column).~Still the time scale gap develops, the fast mode becomes exhausted and $x$ is identified by the CSP pointer as the fast variable.~A valid QSSA can be assumed for $x$.

\begin{figure}[H]
\centering
\includegraphics[width=0.88\linewidth]{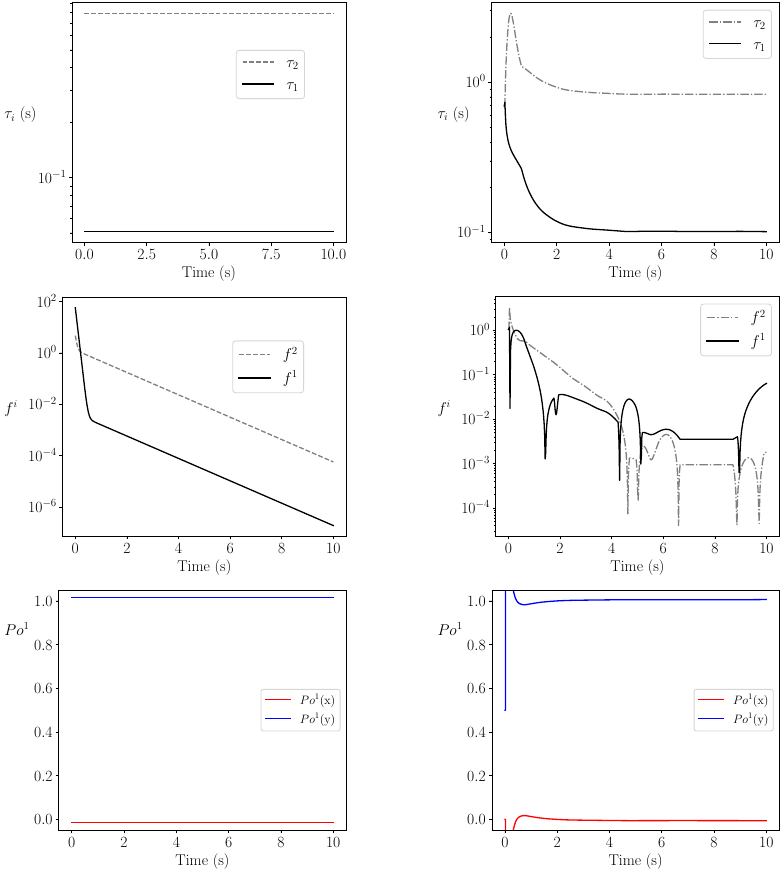}
\caption{
Comparison of the developing timescales (top), CSP amplitudes (middle), and pointers (bottom) for the 2D linear model, as obtained analytically (left) and through the proposed framework (right).~Black line refers to the fast exhausted mode.}
\label{fig:csp_analysis_2Dsnoise}
\end{figure}

\clearpage
\subsection{Michaelis-Menten Model}
\label{app:csp_MM}

\subsubsection{Reference diagnostics}
\label{app:csp_MM_reference_diagnostics}

\begin{figure}[H]
\centering
\includegraphics[width=0.88\linewidth]{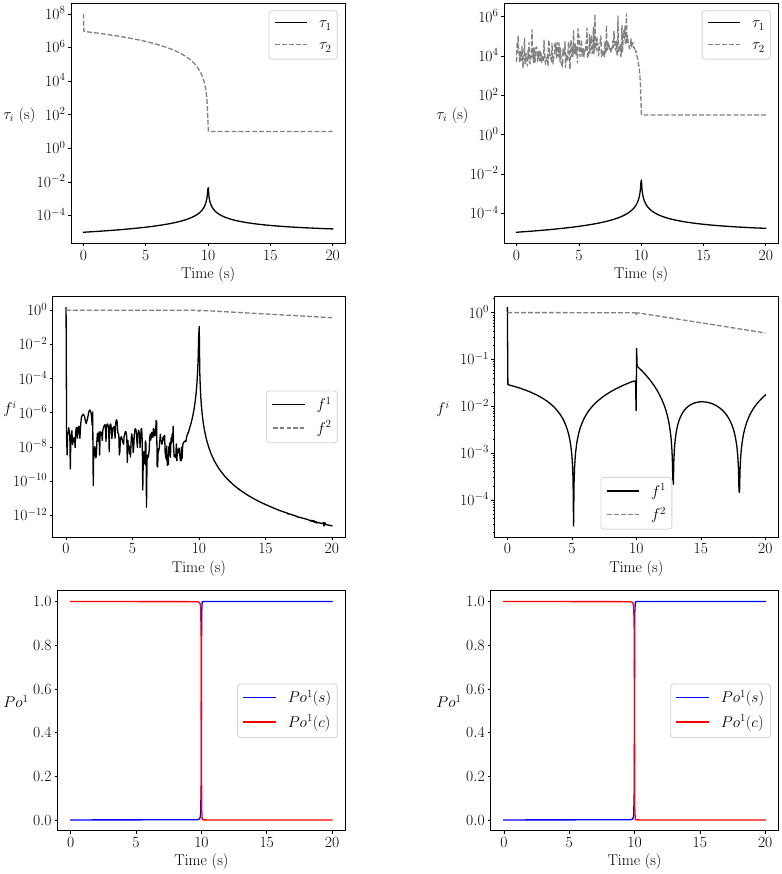}
\caption{
Comparison of the developing timescales (top), CSP amplitudes (middle), and pointers (bottom) for the Michaelis--Menten model, as obtained analytically (left) and through the proposed framework (right). Black line refers to the fast exhausted mode.}
\label{fig:csp_analysis_MM}
\end{figure}

Figure \ref{fig:csp_analysis_MM}~(top) presents the developing time scales of the Michaelis--Menten model as calculated by both the analytical model (left) and the proposed framework (right). A time scale gap is evident in both, indicative of the multi-scale character of the model. It is shown that the amplitude of the fastest dissipative mode becomes negligible (middle), indicating the fast mode exhaustion and forming equilibria, based on which the reduced model can be constructed. Around $t=10s$, a reactivation of the fast amplitude is captured in both setups and a second exhaustion.

The CSP pointer (bottom) was used to identify the variables related the most to the fast dynamics througout the simulations. In both setups, the complex $c$ is identified initially be be the fast variable, indicating thet a valid sQSSA can be applied.~The shift in the dynamics is identified during the reactivation of the fast mode; the substrate $s$ becomes the fast variable in place of the complex $c$. This shift corresponds to the transition between the regions of valid sQSSA and rQSSA.


\subsubsection{Downsampled diagnostics}
\label{app:csp_MMsample}

\begin{figure}[H]
\centering
\includegraphics[width=0.88\linewidth]{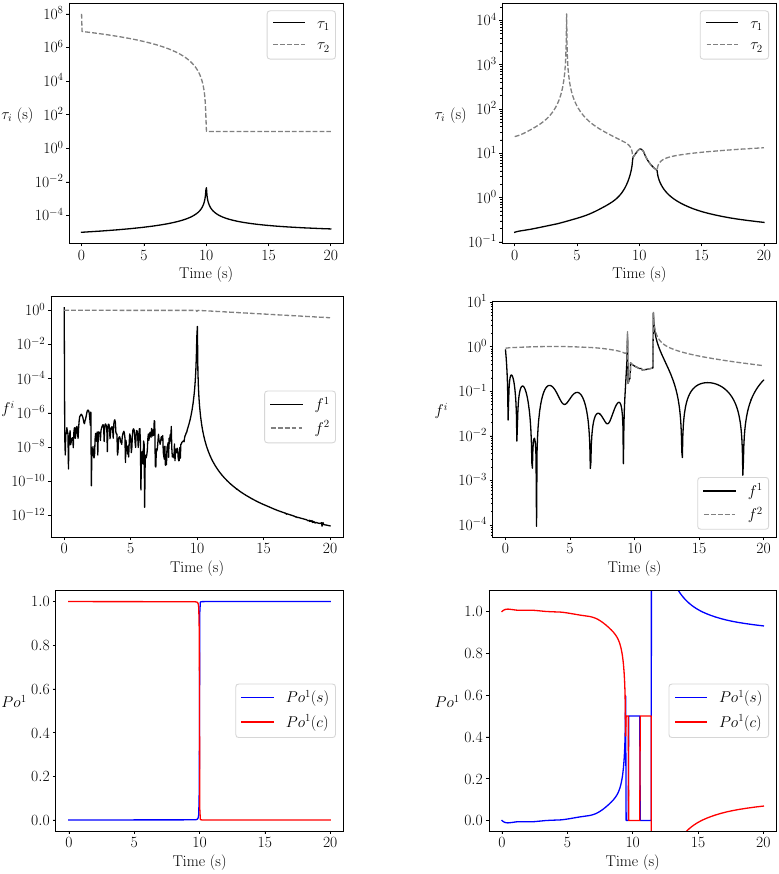}
\caption{
Comparison of the developing timescales (top), CSP amplitudes (middle), and pointers (bottom) for the Michaelis--Menten model, as obtained analytically (left) and through the proposed framework (right). Black line refers to the fast exhausted mode.}
\label{fig:csp_analysis_MMsample}
\end{figure}

The diagnostics obtained on downsampled data by the proposed framework are not affected qualitatively as compared to the ones obtained by the analytical model (left column of Fig.~\ref{fig:csp_analysis_MMsample}).~Still the time scale gap develops, amplitude captures the fast mode that becomes exhausted, reactivated and exhausted again, as a result of the dynamics shift and the CSP pointer identifies first the complex $c$ and then the substrate $s$ and the fast variable in each case, indicating the different valid QSSA that can be assumed.

\subsubsection{Downsampled \& noisy diagnostics}
\label{app:csp_MMnoise}

\begin{figure}[H]
\centering
\includegraphics[width=0.88\linewidth]{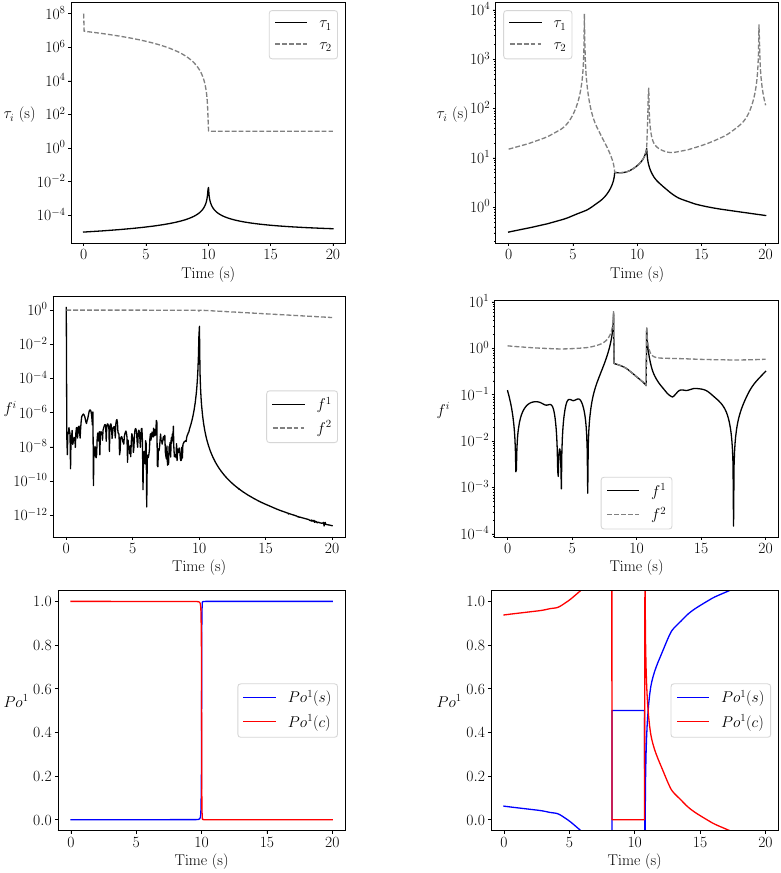}
\caption{
Comparison of the developing timescales (top), CSP amplitudes (middle), and pointers (bottom) for the Michaelis--Menten model, as obtained analytically (left) and through the proposed framework (right). Black line refers to the fast exhausted mode.}
\label{fig:csp_analysis_MMnoise}
\end{figure}

The diagnostics obtained on downsampled and noisy data demonstrate the framework's robustness to noise, showing that key outcomes remain consistent across various perturbation levels for reliable analysis, as compared to the ones obtained by the analytical model (left column).

\subsection{Biofilm Model}
\label{app:csp_biofilm}

\subsubsection{Reference diagnostics}
\label{app:csp_biofilm_reference_diagnostics}

\begin{figure}[H]
\centering
\includegraphics[width=0.88\linewidth]{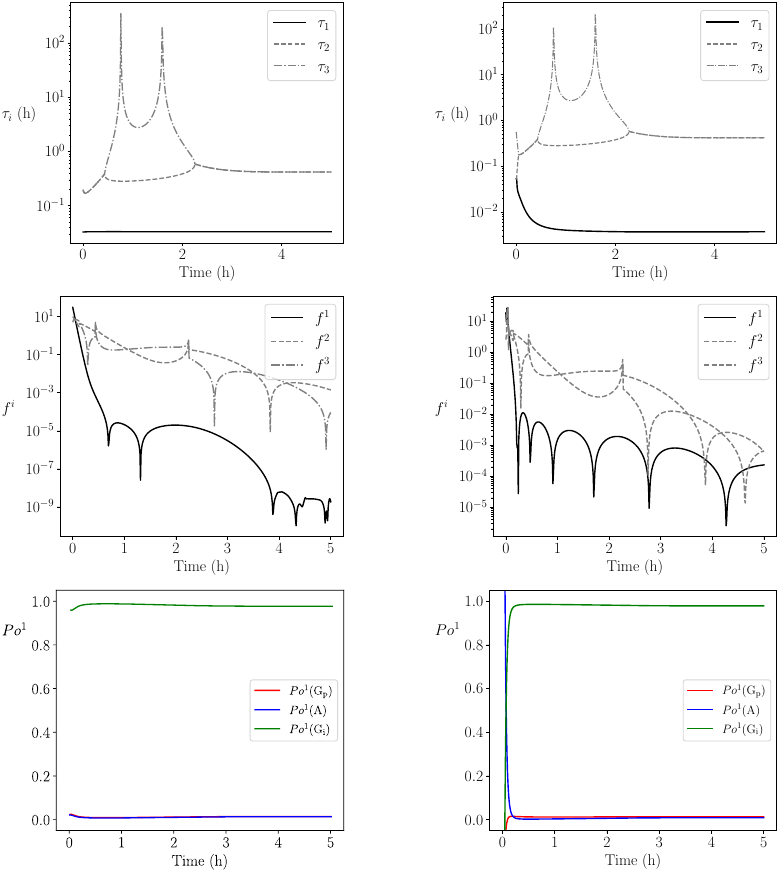}
\caption{Comparison of the developing timescales (top), CSP amplitudes (middle), and pointers (bottom) for the biofilm model, as obtained analytically (left) and through the proposed framework (right). Black line refers to the fast exhausted mode.}
\label{fig:csp_analysis_biofilm}
\end{figure}

The evolution in time of the developing time scales is displayed in Fig. \ref{fig:csp_analysis_biofilm} (top), along with the related modal amplitudes (middle). It is shown that there exists a time scale gap between the fastest and the second faster time scales. CSP modal $Po$ (bottom) identifies that the variable related the most to the fast mode is $G_p$; $G_p$ is the fast variable. Due to the fact that the amplitude of the first mode becomes negligible; i.e. mode 1 becomes exhausted, then a valid QSSA can be assumed for $G_p$.

This QSSA for $G_p$ introduces the description of the SIM, which is derived by considering higher order accuracy as:

\[
G_i^{analytical} \approx G_P \frac{k_4}{k_5} \left[ 1 - \frac{k_1}{k_5} G_E + \frac{k_4}{k_5} + \frac{k_2}{k_5} A \right]
\]

Symbolic regression evaluated over the near-manifold snapshot data produces the following SIM:
\[
G_i^{identified} \approx G_p\Big(\alpha+\beta\,A\Big) + \gamma,
\]
where $(\alpha=0.06353938046,\beta=0.013420349, \gamma =7.412423\times10^{-5})$ are the discovered numerical coefficients.

Substituting the parameter values from Section \ref{app:numerical_settings_Biofilm} the analytical SIM evaluates to::
\[
G_i^{analytical} \approx G_p \left( 0.063498 + 0.011778 A \right),
\] 
which aligns closely with the leading-order structure of the identified expression (omitting the negligible constant offset $\gamma \approx 7.41 \times 10^{-5}$):

\[
G_i^{identified} \approx G_p \left( 0.063539 + 0.013420 A \right).
\]
This strong numerical and structural agreement validates the capability of the proposed data-driven framework.

\clearpage
\subsubsection{Downsampled diagnostics}
\label{app:csp_biofilm_sample}
\begin{figure}[H]
\centering
\includegraphics[width=0.88\linewidth]{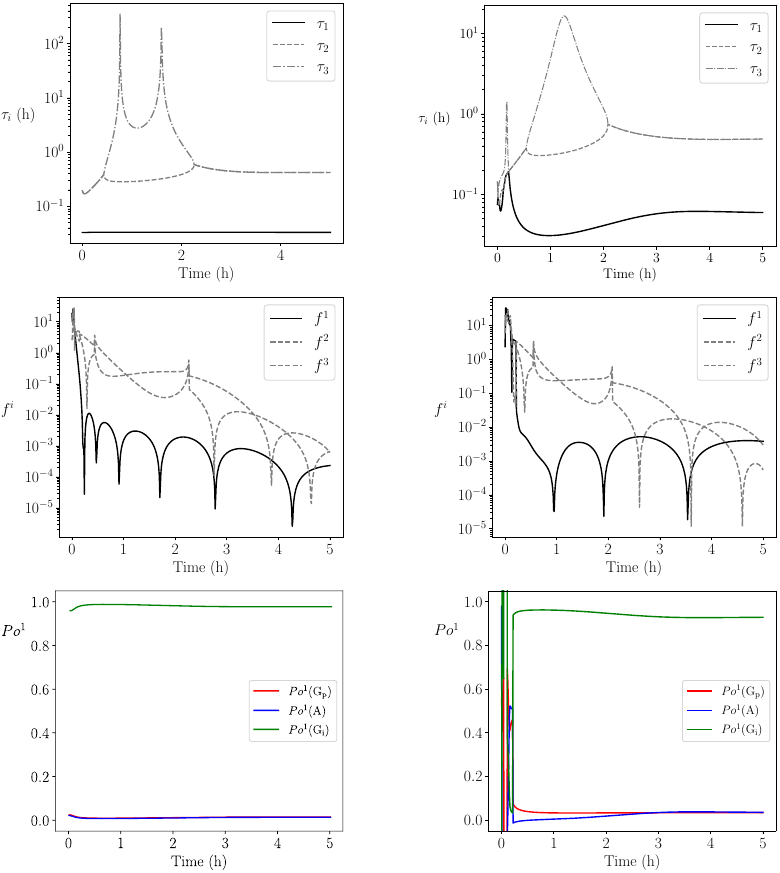}
\caption{
Comparison of the developing timescales (top), CSP amplitudes (middle), and pointers (bottom) for the biofilm model, as obtained analytically (left) and through the proposed framework (right). Black line refers to the fast exhausted mode.}
\label{fig:csp_analysis_biofilm_sample}
\end{figure}

The diagnostics obtained on downsampled data by the proposed framework are not affected qualitatively as compared to the ones obtained by the analytical model (left column of Fig.~\ref{fig:csp_analysis_biofilm_sample}).~Still the time scale gap develops, modal ampliitude captures the fast exhausted mode and the CSP pointer identifies $G_p$ as the fast variable; a valid QSSA can be assumed for $G_p$.

\clearpage
\subsubsection{Downsampled \& noisy diagnostics}
\label{app:csp_biofilm_noise}

\begin{figure}[H]
\centering
\includegraphics[width=0.88\linewidth]{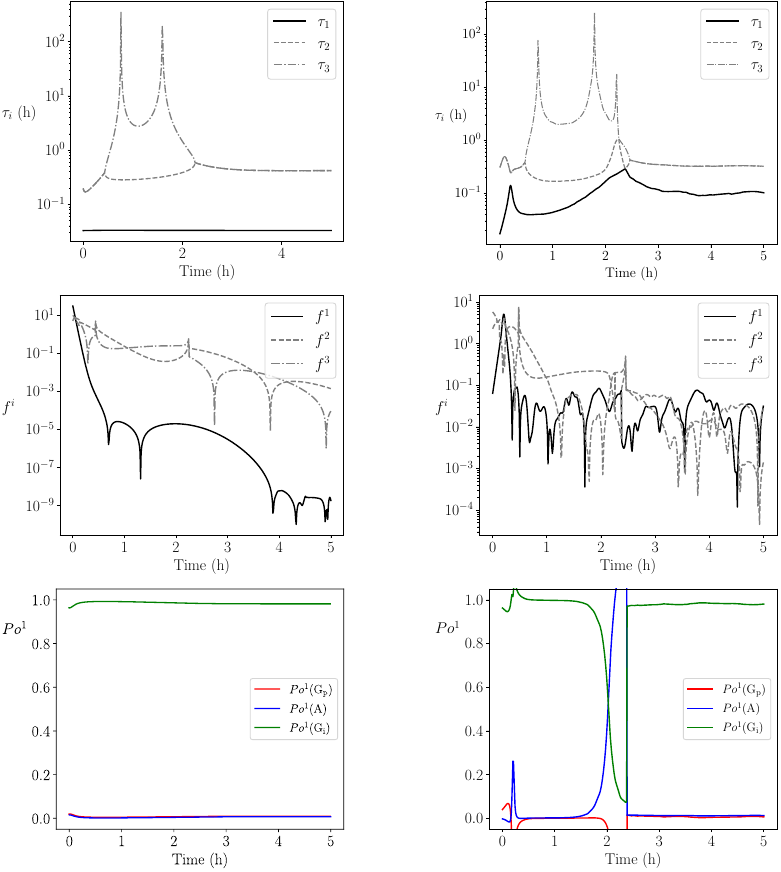}
\caption{
Comparison of the developing timescales (top), CSP amplitudes (middle), and pointers (bottom) for the biofilm model, as obtained analytically (left) and through the proposed framework (right). Black line refers to the fast exhausted mode.}
\label{fig:csp_analysis_biofilm_noise}
\end{figure}

The diagnostics obtained on downsampled and noisy data by the proposed framework are affected qualitatively as compared to the ones obtained by the analytical non-noisy model (left column of Fig.~\ref{fig:csp_analysis_biofilm_noise}).~The time scale gap develops, but to a smaller degree, the modal ampliitude is affected by the noise in distingiuishing the exhausted mode. Hoever, the CSP pointer identifies $G_p$ to be mostly associated witht the fast mode; $G_p$ is identified as the fast variable and a valid QSSA can be assumed for $G_p$.

\subsection{Target-mediated drug disposition model}
\label{app:csp_TMDD}

\subsubsection{Reference diagnostics}
\label{app:csp_TMDD_reference_diagnostics}

\begin{figure}[H]
\centering
\includegraphics[width=0.88\linewidth]{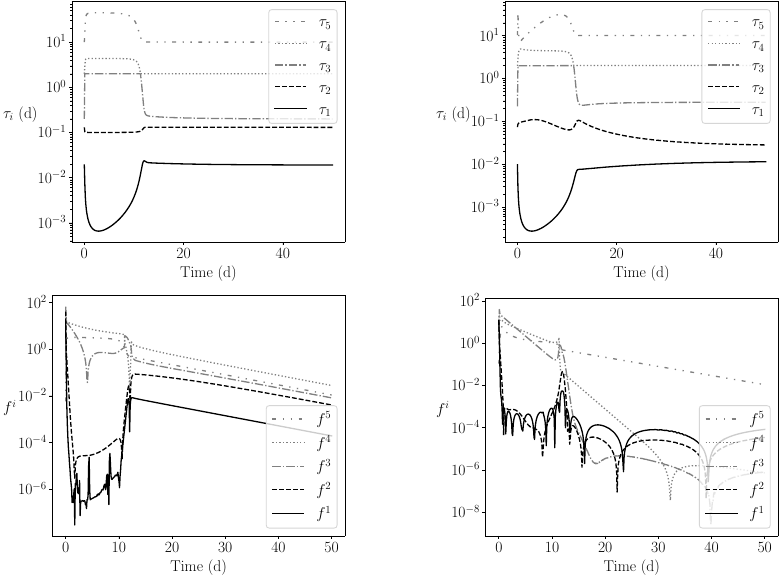}
\caption{Comparison of the developing timescales (top) and CSP modal amplitudes (bottom) for the TMDD model, as obtained analytically (left) and through the proposed framework (right). Black lines refer to the fast exhausted modes.}
\label{fig:csp_analysis_TMDD}
\end{figure}

\begin{figure}[H]
\centering
\includegraphics[width=0.88\linewidth]{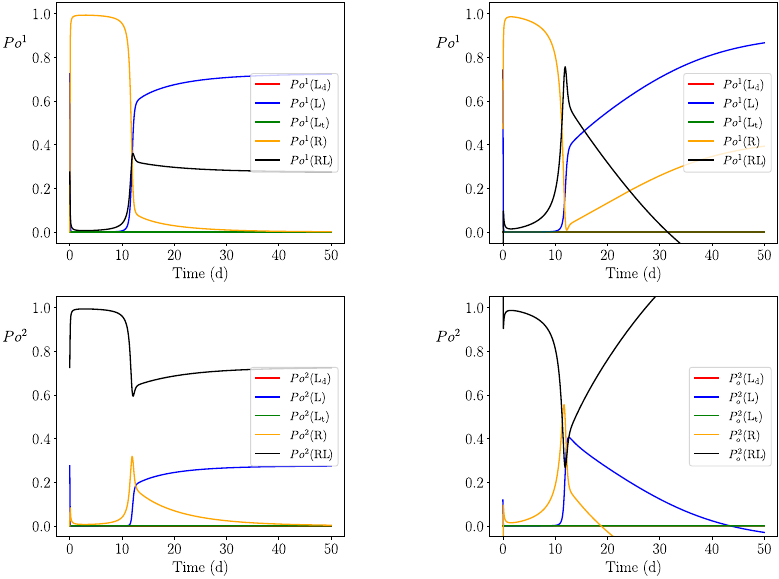}

\caption{
Modal pointer associated with the first (top) and second (bottom) CSP mode for the TMDD model, obtained analytically (left) and through the proposed framework (right).}
\label{fig:csp_analysis_TMDD_b}
\end{figure}

Figure \ref{fig:csp_analysis_TMDD} (top) presents the developing time scales of the target-mediated drug disposition (TMDD) model as calculated by both the analytical model (left) and the proposed framework (right). A clear time scale gap is evident in both, indicative of the multiscale character of the model. During the initial period (from $t = 0$ to $10\text{ d}$), the modal amplitudes of the first two modes become negligible (middle), indicating fast mode exhaustion and the formation of equilibria, based on which a reduced model can be constructed. The CSP pointers (Fig. \ref{fig:csp_analysis_TMDD_b}) identify that variables $R$ and $RL$ are directly associated with these first two exhausted modes; $R$ and $RL$ are the fast variables for which a valid QSSA can be assumed during this initial phase.

After $t = 10\text{ d}$, the diagnostics show that a time scale gap exists between the third and fourth time scales. Correspondingly, the modal amplitudes of the first three modes become exhausted, expanding the equilibrium structure for a subsequent reduced model. The CSP pointer identifies three fast variables, namely $R$, $RL$, and $L$, for which a valid QSSA can be assumed in this subsequent regime.

\vfill
\clearpage

\label{LastPageSupp3}

\end{document}